\documentclass[english, aps, prl, reprint, notitlepage, superscriptaddress, longbibliography]{revtex4-1}
\usepackage[T1]{fontenc}
\usepackage[latin9]{inputenc}
\usepackage{amsbsy}
\usepackage{amstext}
\usepackage{graphicx}
\usepackage{esint}
\usepackage{babel}
\usepackage{xcolor}
\usepackage[normalem]{ulem}
\usepackage{url}

\usepackage[breaklinks=true]{hyperref}
\usepackage{breakcites}
\usepackage{breakurl}

\begin{document}

\title{Bilevel optimization in flow networks -- A message-passing approach}

\author{Bo Li}
\affiliation{Non-linearity and Complexity Research Group, Aston University, Birmingham,
B4 7ET, United Kingdom}
\affiliation{School of Science, Harbin Institute of Technology (Shenzhen), Shenzhen, 518055, China}

\author{David Saad}
\affiliation{Non-linearity and Complexity Research Group, Aston University, Birmingham,
B4 7ET, United Kingdom}

\author{Chi Ho Yeung}
\affiliation{Department of Science and Environmental Studies, The Education University of Hong Kong, 10 Lo Ping Road, Taipo, Hong Kong}

\begin{abstract}
Optimizing embedded systems, where the optimization of one depends on the state of another, is a formidable computational and algorithmic challenge, that is ubiquitous in real world systems. We study flow networks, where bilevel optimization is relevant to traffic planning, network control and design, and where  flows are governed by an optimization requirement subject to the network parameters. We employ message-passing algorithms in  flow networks with sparsely coupled structures to adapt network parameters that govern the network flows, in order to optimize a global objective. We demonstrate the effectiveness and efficiency of the approach on randomly generated graphs.
\end{abstract}

\maketitle

Many problems in science and engineering involve hierarchical optimization, whereby some of the variables cannot be freely varied but are governed by another optimization problem~\cite{Sinha2018}.
As a motivating example, consider the task of designing a network (e.g., a road or communication network) that maximizes the throughput of commodities or information flow. While the designer controls the network parameters (upper-level optimization), traffic flows are determined by the network users who maximize their own benefit (lower-level optimization)~\cite{Wardrop1952}.
Therefore, the designer needs to adapt the network intricately, taking into account the reaction of network users. 
Similarly, many physical systems admit a certain extremization principle for given controllable system parameters, e.g., minimal free energy in thermal equilibrium~\cite{Plischke2006}, electric flows in resistor networks that minimize dissipation~\cite{Thomson2009,Doyle1984} and, entropy maximization and parameter optimization that are used across disciplines in inference and learning tasks~\cite{Janes1957,Murphy2012}.
Adapting system parameters to extremize a given objective requires bilevel optimization, which considers both system parameters and the inherent optimization of the physical variables.

Bilevel optimization is intrinsically difficult to solve~\cite{Colson2007}. In fact, even the simple instance where both levels are linear programming tasks is NP-hard~\cite{Jeroslow1985,Hansen1992}. Generic methods for bilevel optimization include (i) bilevel programming approach, by expressing the lower-level optimization problem as nonlinear constraints and solving the bilevel problem as global optimization~\cite{Bard1982,Bard1988}; (ii) gradient-descent method by computing the descent direction of the upper-level objectives while keeping the valid lower-level state variables~\cite{Kolstad1990,Savard1994}.
The former introduces complicated nonlinear constraints, making the reduced single-level problem difficult in general, while the latter can be challenging in computing the descent direction~\cite{Colson2007}. Moreover, such generic methods do not utilize existing system structure to simplify the task.

In this Letter, we develop message-passing (MP) algorithms to tackle bilevel optimization in sparse flow networks. The  advances presented in this work are three-fold: (i) the derived MP algorithms are intrinsically distributed, scalable, and generally efficient; (ii) they are applicable to bilevel optimization problems with combinatorial constraints, which are difficult for generic bilevel programming approaches; (iii) these algorithms can successfully deal with nonsmooth flow problems, having  potential applications  for transport based approaches in machine learning~\cite{Rustamov2018, Essid2018, Peyre2019}.

\textit{Routing Game.} We focus on a network planning problem in the routing game setting, widely used in modeling route choices of drivers~\cite{Patriksson2015}. Users on the road network make their route choices in a selfish and rational manner, where the corresponding Nash equilibrium is generally not the most beneficial for the global utility, measured by the total travel time of all users~\cite{Wardrop1952,Roughgarden2005book}.
The operator's task is to set the appropriate tolls or rewards on network edges to reduce the total travel time while taking into account the reactions
of users to the tolls~\cite{Beckmann1956,Smith1979,Cole2006}.
Recently, the idea of reducing traffic congestion by economic incentives to influence drivers' behaviors has regained interest~\cite{Colak2016,ERP_Singapore,Barak2019}, partly due to the deployment of smart devices and data availability~\cite{AlonsoMora2017,Lim2011,BoLiPRR2020}. Here, we focus on the algorithmic aspect of toll optimization.

The road network is represented by a directed graph $G(V,E)$, where $V$ is the set of nodes (junctions) and $E$ the set of directed edges (unidirectional roadways), having one connected component. 
Users routing from an origin node $i_{0}$ to a destination node $\mathcal{D}$ would select a path $\mathcal{P}=((i_{0},i_{1}),(i_{1},i_{2}),...,(i_{n-2},i_{n-1}),(i_{n-1},\mathcal{D}))$ by minimizing their total travel time $\sum_{e\in\mathcal{P}}\ell_{e}(x_{e})$, or alternative cost, where the edge flow $x_{e}$ represents the number of users choosing edge $e$ and $\ell_{e}(x_{e})$ is the corresponding latency function. It is assumed that $\ell_{e}$ is monotonically increasing with the edge flow $x_{e}$. 
The social cost is defined as the total travel time of all users $H=\sum_{e\in E}x_{e}\ell_{e}(x_{e})$, which is the overall objective of the bilevel optimization problem.

We consider the limit of a large number of users, where each user controls an infinitesimal fraction of the overall traffic, such that the edge flow $x_{e}$ is a continuous variable. This is termed the nonatomic game setting~\cite{Roughgarden2005book}.
As the equilibrium reached by the selfish decisions of users does not generally achieve the lowest social cost, we seek to place tolls $\{\tau_{e}\}$ on edges to influence users' route choices. Gauging the monetary penalty at the same scale as latency, users will choose a path $\mathcal{P}$ that minimizes the combined total journey cost in latency and tolls $\sum_{e\in\mathcal{P}}\big[\ell_{e}(x_{e})+\tau_{e}\big]$. 
If tolls can be placed freely on all edges, marginal cost pricing is known to induce socially optimal flow for nonatomic games~\cite{Smith1979}. However, it is usually infeasible to set an unbounded toll on every road, which renders marginal cost pricing less applicable. We therefore consider restricted tolls $0\leq\tau_{e}\leq  \tau_{e}^{\max}$; an edge $e$ is not chargeable when $\tau_{e}^{\max}=0$. For simplicity, we do not consider the income from tolls to contribute to the social cost~\cite{Karakostas2005}.
In total, $\Lambda_{i}$ users are traveling from node $i$ to a universal destination $\mathcal{D}$, where the case with multiple destinations is discussed in the supplemental material (SM)~\cite{BoLi2021sup}. 
The resulting edge flows satisfy the non-negativity $x_{e}\geq0$ and the flow conservation constraints
\begin{equation}
R_{i}=\Lambda_{i}+\sum_{e\in\partial_{i}^{\text{in}}}x_{e}-\sum_{e\in\partial_{i}^{\text{out}}}x_{e}=0,\label{eq:flow_conservation}
\end{equation}
where $\partial_{i}^{\text{in}}$ and $\partial_{i}^{\text{out}}$ are the sets of incoming and outgoing edges adjacent to node $i$.
It has been established that the edge flows in user equilibrium (i.e., the Wardrop's equilibrium~\cite{Wardrop1952}) can
be obtained by minimizing a potential function $\Phi = \sum_{e\in E}\phi_{e}(x_{e}):=\sum_{e\in E}\int_{0}^{x_{e}}  [\ell_{e}(y)+\tau_{e}] \mathrm{d}y$ subject to the constraints of Eq.~(\ref{eq:flow_conservation}) ~\cite{Monderer1996,BarGera2002}.
We emphasize that the potential function $\Phi(\boldsymbol{x})$ only plays an auxiliary role in defining the equilibrium flows; the values of $\Phi$ do \emph{not} correspond to the routing costs of users.

\begin{figure}
\hspace*{-2em}\includegraphics[scale=1.1]{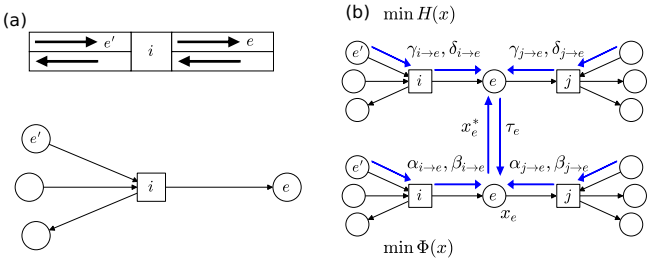}
\caption{(a) Top: a directed road network section with a junction node $i$. Bottom: the corresponding factor graph representation; node $i$ is a factor node and is marked by a square. (b) Bilevel MP for toll planning. Blue arrows indicate the directions of messages.
The equilibrium flow $x_{e}^{*}$ is determined in the lower level, while the toll $\tau_{e}$ is set in the upper level.
\label{fig:routing_game_bilevel_illustrate}}
\end{figure}

The lower-level optimization is a nonlinear min-cost flow problem, where edge flows are coupled through the conservation constraints in Eq.~(\ref{eq:flow_conservation}), represented as factor nodes in Fig.~\ref{fig:routing_game_bilevel_illustrate}(a).
We employ the MP approach developed in Ref.~\cite{Wong2007}
to tackle the nonlinear optimization problem. It turns the global optimization of the potential into a local computation of the following message functions
\begin{equation}
\Phi_{i\to e}(x_{e})=\min_{\{x_{e'}\geq0\}|R_{i}=0}\sum_{e'\in\partial i\backslash e}\bigg[\Phi_{e'\to i}(x_{e'})+\phi_{e'}(x_{e'})\bigg],\label{eq:cavity_Phi_ie}
\end{equation}
where $\partial i=\partial_{i}^{\text{in}}\cup\partial_{i}^{\text{out}}$
and $\Phi_{i\to e}(x_{e})$ relates to the optimal potential function contributed by the flows adjacent to node $i$ where the flow on edge $e$ is set to $x_{e}$, taking into account flow conservation at node $i$. 
In Eq.~(\ref{eq:cavity_Phi_ie}), denoting $e'=(k,i)$, we can write $\Phi_{e'\to i}(x_{e'})=\Phi_{k\to e'}(x_{e'})$; therefore only factor-to-variable messages are needed. The message $\Phi_{k\to e'}(x_{e'})$ can be obtained recursively by an expression similar to Eq.~(\ref{eq:cavity_Phi_ie}), but using the incoming messages from its upstream edges $\{l\to k|(l,k)\in\partial k\backslash i\}$.
Upon computing the messages iteratively until convergence, we can determine the equilibrium flow $x_{e}^{*}$ on edge $e=(i,j)$ by minimizing the edgewise full energy dictated by the nonlinear cost $\phi_{e}(x_{e})$ and messages from both ends of edge $e$, defined as $\Phi_{e}^{\text{full}}(x_{e})=\Phi_{i\to e}(x_{e})+\Phi_{j\to e}(x_{e})+\phi_{e}(x_{e})$.

\begin{figure}
\hspace*{-0.5em}\includegraphics[scale=1.1]{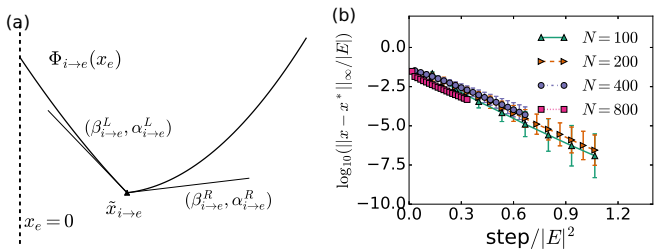}
\caption{(a) A nonsmooth message function $\Phi_{i \to e}(x_{e})$ with one breakpoint. (b) Convergence of the single-level MP algorithm for computing equilibrium flows in routing games in random regular graphs with degree 3 of different sizes $N=|V|$. An affine latency model $\ell_{e}(x_{e})=t_{e}(1+sx_{e}/c_{e})$ is considered, where $t_{e}$ and $c_{e}$ are the free traveling time and edge capacity, respectively, while $s$ is a sensitivity measure of latency to congestion~\cite{Roughgarden2005book}. Random sequential schedule of MP updates has been used.
}
\label{fig:routing_single_MP}
\end{figure}

This algorithm can be demanding when different values of $x_{e}$ are needed to determine the profile of the message $\Phi_{i\to e}(x_{e})$.
To reduce the computational cost, we consider the approximation of the message in the vicinity of some working point $\tilde{x}_{i\to e}$ as 
\begin{equation}
\Phi_{i\to e}(\tilde{x}_{i\to e}+\varepsilon_{e})\approx\text{\ensuremath{\Phi_{i\to e}}(\ensuremath{\tilde{x}_{i\to e}}})+\beta_{i\to e}\varepsilon_{e}+\frac{1}{2}\alpha_{i\to e}\big(\varepsilon_{e}\big)^{2},\label{eq:Phi_ie_approx}
\end{equation}
where $\beta_{i\to e}$ and $\alpha_{i\to e}$ are the first and second derivatives of $\Phi_{i\to e}$ evaluated at $\tilde{x}_{i\to e}$, assuming the derivatives exist. 
For a particular $\tilde{x}_{i\to e}$, the computation of the message function $\Phi_{i\to e}(x_{e})$ in Eq.~(\ref{eq:cavity_Phi_ie}) reduces to the optimization of $\beta_{i\to e}$ and $\alpha_{i\to e}$ by using $\{\tilde{x}_{k\to e'},\beta_{k\to e'},\alpha_{k\to e'}|e'=(k,i)\in\partial i\backslash e\}$. 
The working point $\tilde{x}_{i\to e}$ is updated by pushing it towards the minimizer $x_{e}^{*}$ of the full energy $\Phi_{e}^{\text{full}}(x_{e})$ gradually~\cite{BoLi2021sup}.
The iterative updates of the coefficients $\{ \beta_{i \to e}, \alpha_{i \to e} \}$ and the working points $\{ \tilde{x}_{i \to e} \}$ constitute a perturbative version of the original MP algorithm, which only requires to keep track of a few coefficients rather than the full profile of $\Phi_{i\to e}$, making it tractable~\cite{Wong2007}. It has been shown to work remarkably well in many network flow problems~\cite{Wong2016}, while the algorithm may not converge in problems with nonsmooth characteristics~\cite{Wong2007}. 
We discover that the non-negativity constraints on flows can result in a nonsmooth message function $\Phi_{i\to e}(x_{e})$, which makes the approximation of Eq.~(\ref{eq:Phi_ie_approx}) inadequate. 
One solution is to approximate $\Phi_{i\to e}(x_{e})$ by a continuous and piecewise quadratic function with at most two branches, where each branch $m$ is a quadratic function governed by three coefficients $\{ \tilde{x}_{i\to e}, \beta_{i\to e}^{(m)}, \alpha_{i\to e}^{(m)}\}$, as illustrated in Fig.~\ref{fig:routing_single_MP}(a).
Taking into account the nonsmooth structures, MP algorithms converge well even in loopy networks and provide the correct solutions~\cite{BoLi2021sup}.
We demonstrate the case of random regular graphs (RRG) with degree 3 in Fig.~\ref{fig:routing_single_MP}(b).

For bilevel optimization, we notice that the cost function of the upper layer $H(x)$ has a similar structure as $\Phi(x)$. Therefore, one can apply a similar MP procedure as $H_{i\to e}(x_{e})=\min_{\{x_{e'}\}|R_{i}=0}\sum_{e'\in\partial i\backslash e}[H_{e'\to i}(x_{e'})+x_{e'}\ell_{e}(x_{e'})]$.
The message $H_{i\to e}(x_{e})$ can also be approximated by a piecewise quadratic function with at most one break point, where each branch $m$ has the form $H_{i\to e}^{(m)}(\tilde{x}_{i\to e}+\varepsilon_{e})\approx\text{\ensuremath{H_{i\to e}^{(m)}}(\ensuremath{\tilde{x}_{i\to e}}})+\gamma_{i\to e}^{(m)}\varepsilon_{e}+\frac{1}{2}\delta_{i\to e}^{(m)}\big(\varepsilon_{e}\big)^{2}$.
As the equilibrium state is determined in the lower level, the working points $\{\tilde{x}_{i\to e}\}$ in the lower-level MP are also used for the upper level. 
The landscape of the edgewise full cost $H_{e}^{\text{full}}(x_{e})=H_{i\to e}(x_{e})+H_{j\to e}(x_{e})+x_{e}\ell(x_{e})$ provides the information for setting the toll. 
Specifically, the toll is updated by $\min_{\tau_{e}}H_{e}^{\text{full}}(x_{e}^{*}(\tau_{e}))$, where the toll-dependent equilibrium flow $x_{e}^{*}$ is provided by the lower-level messages.
In practice, an approximate $H_{e}^{\text{full}}$ is sufficiently informative for updating tolls.
The basic structure of such bilevel MP is illustrated in Fig.~\ref{fig:routing_game_bilevel_illustrate}(b), while details are provided in the SM~\cite{BoLi2021sup}.

\begin{figure}
\hspace*{-0.5em}\includegraphics[scale=1.05]{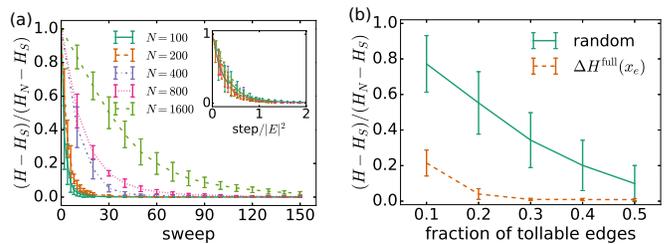}
\caption{Bilevel MP algorithm for routing games on RRG. (a) Effect of tolls on the fractional social cost reduction $(H(\boldsymbol{x}^{*}(\boldsymbol{\tau}))-H_{S})/(H_{N}-H_{S})$, where $H_{S}$ and $H_{N}$ represent the social costs at the social optimum and the Nash equilibrium without tolls.
Tolls $\boldsymbol{\tau}$ are recorded during bilevel MP updates, based on the resulting equilibrium flows $\boldsymbol{x}^{*}(\boldsymbol{\tau})$ and social cost $H(\boldsymbol{x}^{*}(\boldsymbol{\tau}))$. 
Each data point is the average of 10 different problem realizations. Each sweep consists of $40|E|$ local MP steps and $100$ edgewise toll updates in a random sequential schedule. A fixed number of sweeps without toll updates are performed to warm up the system. Inset: panel (a) with $x$-axis as MP steps rescaled by $|E|^{2}$. 
(b) Fractional cost reduction as a function of the fraction  of tollable edges on an RRG with $N=200$. A random selection of edges to be charged is compared with selections based on edgewise full cost reduction $H_{e}^{\text{full}}(x_{e}^{*})$.
}
\label{fig:routing_bilevel_results}
\end{figure}

We demonstrate the effectiveness of the proposed bilevel MP algorithm for tasks on RRG in Fig.~\ref{fig:routing_bilevel_results}(a), where the set-up is the same as in Fig.~\ref{fig:routing_single_MP}(b).
Experiments on other networks and the cases of multiple destinations are discussed in the SM~\cite{BoLi2021sup}.
Although bilevel message-passing does not generally converge to a set of \emph{unique} optimal tolls due to the non-convex nature of the problem, we found that the social costs are reduced when tolls are updated during MP. The scaling relation in the inset of Fig.~\ref{fig:routing_bilevel_results}(a) empirically indicates that the number of updates is $O(|E|^{2})$ for achieving a given cost reduction.
Moreover, the MP algorithm can be implemented in a fully distributed manner, unlike the generic global optimization approach~\cite{BoLi2021sup}. 
Note that we have utilized the special set-up of routing games here, where the social optimum $H_{S} = \min_{x} H(x)$ can be obtained \textit{a priori} for this benchmark. Such information may be unavailable in other bilevel-optimization problems.
The toll optimization problem can also be tackled by the bilevel programming approach~\cite{Bard1982,Bard1988}; however, it requires a treatment with mixed integer programming, which is centralized and generally not scalable, unlike the MP approach~\cite{BoLi2021sup}.

\textit{Combinatorial Problems.} In practice, it may be infeasible to charge for every edge, but desirable to choose a subset of tollable edges for toll-setting~\cite{Hoefer2008}, which is a difficult combinatorial optimization problem. 
As the cost landscape is manifested locally by the message functions, we heuristically select the tollable edges according to the largest possible reduction in edgewise full cost $H_{e}^{\text{full}}(x_{e}^{*})$ due to tolling, which effectively selects the chargeable links as seen in Fig.~\ref{fig:routing_bilevel_results}(b). 
Such combinatorial problems are generally very difficult for traditional bilevel-optimization methods, while MP algorithms can provide approximate solutions in some scenarios.

Another important class of combinatorial problems is the atomic games which consider integer flow variables $\{ x_{e} \}$~\cite{Rosenthal1973}. In principle, atomic games can be solved via the same MP procedure as in Eq.~(\ref{eq:cavity_Phi_ie}), where the message $\Phi_{i \to e}(x_{e})$ is defined on a one dimensional grid. Using the techniques in Refs.~\cite{Yeung2012,Yeung2013IEEE,Bacco2014,Yeung2013,Po2021}, the MP approach provides a scalable algorithm to approximately tackle the difficult combinatorial optimization of atomic games in a single level; it can also solve instances of the bilevel toll-optimization problems. However, its performance is sub-optimal in large networks and for cases with heavy loads~\cite{BoLi2021sup}.
Nevertheless, we found some interesting patterns of the optimal tolls in a realistic test case network using this method~\cite{BoLi2021sup}.

\textit{Flow Control.} 
We consider the problem of tuning network flows to achieve certain functionality. 
In this example, resources need to be transported from source nodes to destination  along edges in an undirected network $G(V,E)$, where the equilibrium flows $\{x_{ij}^{*}\}$ minimize the transportation cost $C=\sum_{(i,j)\in E}\frac{1}{2}r_{ij}x_{ij}^{2}$, subject to flow conservation constraints similar to Eq.~(\ref{eq:flow_conservation}).
The major difference of this model from routing games is that the network is undirected, where edge $(i,j)$ can accommodate either the flow from node $j$ to $i$ or $i$ to $j$.
The objective is to control the parameters $\{r_{ij}\}$ to reduce or increase the flows on some edges.
The task of reducing edge flows has applications in power grid congestion mitigation in the direct-current (DC) approximation~\cite{Wood2013}, where $r_{ij}$ is related to the reactance of edge $(i,j)$, controllable through devices in a flexible alternating current transmission system (FACTS)~\cite{Zhang2006}. On the other hand, the task of increasing certain edge flows has been used to model the tunability of network functions, which is applicable in mechanical and biological networks~\cite{Rocks2019} as well as learning machines in metamaterials~\cite{Stern2021}. 

\begin{figure}
\includegraphics[scale=1.05]{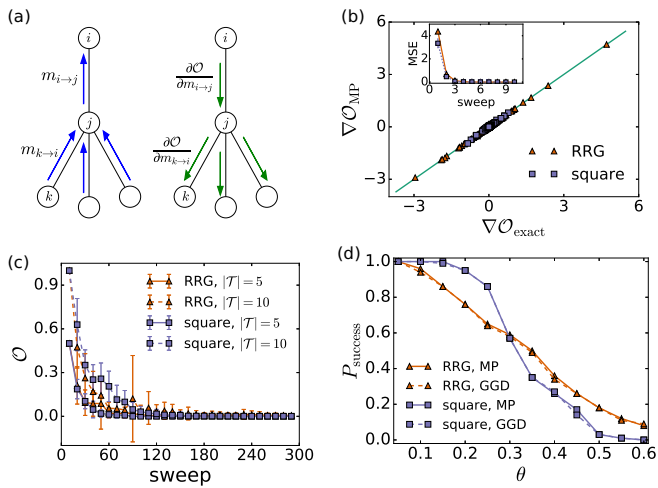}
\caption{Bilevel optimization for flow control. An RRG ($N=200$, degree $3$) and a square lattice of size $15\times15$ are considered. The source and destination nodes, and the targeted edges are randomly selected. (a) Left: MP for solving the lower-level equilibrium flow problem. Right: Computing gradients of the upper-level objective function $\mathcal{O}$.
(b) Comparison of the gradients at initial $\boldsymbol{r}$ computed by the MP approach (obtained by fixing $\boldsymbol{r}$ and passing messages $\{m_{i\to j}\}$ and gradients $\{\frac{\partial\mathcal{O}}{\partial m_{i\to j}}\}$) and the GGD approach, with $|\mathcal{T}|=5,\theta=0.1$. Inset: mean square error (MSE) of the gradients by the MP approach during iterations, in comparison to the GGD approach. 
Each sweep consists of $4|E|$ local MP steps. (c) MP for minimizing the upper-level objective function $\mathcal{O}$ with $\theta=0.1$, where one randomly selected control parameter is updated following the descent direction every $4|E|/10$ steps. (d) Fraction of successfully tuned cases (satisfying $\mathcal{O}=0$) $P_{\text{success}}$ out of $100$ different problem realizations of source/destination nodes, with $|\mathcal{T}|=5$, as a function of the threshold $\theta$. 
\label{fig:bilevel_flow_control_results}}
\end{figure}

As an example, we consider the task of flow control such that the relative increments of flows on the targeted edges $\mathcal{T}$ exceed a limit $\theta$~\cite{Rocks2019}, i.e., $\rho_{ij}=\frac{|x_{ij}|-|x_{ij}^{0}|}{|x_{ij}^{0}|}-\theta\geq0,\forall(i,j)\in\mathcal{T}$ (with $x_{ij}^{0}$ being the flow prior to tuning). 
It can be achieved by minimizing the hinge loss (upper-level objective) $\mathcal{O}=\sum_{(i,j)\in\mathcal{T}}-\rho_{ij}\Theta(-\rho_{ij}) =: \sum_{(i,j)\in\mathcal{T}} \mathcal{O}_{ij} $, where $\Theta(\cdot)$ is the Heaviside step function. 
The task of congestion mitigation in power grids can be studied similarly. We adopt the usual MP algorithm to compute the equilibrium flows as
\begin{equation}
C_{i\to j}(x_{ij})=\min_{\{x_{ki}\}|R_{i}=0}\bigg[\frac{1}{2}r_{ij}x_{ij}^{2}+\sum_{k\in \mathcal{N}_{i}\backslash j}C_{k\to i}(x_{ki})\bigg],\label{eq:cavity_C_ij}
\end{equation}
where $\mathcal{N}_{i}$ is the set of neighboring nodes adjacent to node $i$. The definition of the message $C_{i\to j}(x_{ij})$ differs from the one of Eq.~(\ref{eq:cavity_Phi_ie}) in that it includes
the interaction term on edge $(i,j)$, which yields a more concise update rule here. Similar to Eq.~(\ref{eq:Phi_ie_approx}), we approximate the message function by a quadratic form $C_{i\to j}(x_{ij})=\frac{1}{2}\alpha_{i\to j}(x_{ij}-\hat{x}_{i\to j})^{2}+\text{const}$, such that the optimization in Eq.~(\ref{eq:cavity_C_ij}) reduces to the computation of the real-valued messages $m_{i\to j}\in\{\alpha_{i\to j},\hat{x}_{i\to j}\}$ by passing the upstream messages $\{m_{k\to i}\}_{k\in\mathcal{N}_{i}\backslash j}$, as illustrated on the left panel of Fig.~\ref{fig:bilevel_flow_control_results}(a)~\cite{BoLi2021sup}. Upon convergence, the equilibrium flow $x_{ij}^{*}$ can be obtained by minimizing the edgewise full cost $C_{ij}^{\text{full}}(x_{ij})=C_{i\to j}(x_{ij})+C_{j\to i}(x_{ij})-\frac{1}{2}r_{ij}x_{ij}^{2}$.

The variation of the control parameters $\{r_{ij}\}$ will impact on the messages $\{m_{i\to j}\}$, which in turn affects the equilibrium flows $\boldsymbol{x}^{*}$ and therefore the upper-level objective $\mathcal{O}(\boldsymbol{x}^{*})$. Specifically, one considers the effect of the change of $r_{ij}$ on the targeted edge flows $\{x_{pq}^{*}\}_{(p,q)\in\mathcal{T}}$, derived by computing the gradient $\frac{\partial\mathcal{O}}{\partial m_{i\to j}}$.
The targeted edges provide the boundary conditions as $\frac{\partial\mathcal{O}_{pq}}{\partial m_{p\to q}}=\frac{\partial\mathcal{O}_{pq}}{\partial x_{pq}^{*}}\frac{\partial x_{pq}^{*}}{\partial m_{p\to q}},\forall(p,q)\in\mathcal{T}$.
As the messages from node $i$ to $j$ are functions of the upstream messages, i.e., $m_{i\to j}=m_{i\to j}(\{m_{k\to i}\}_{k\in\mathcal{N}_{i} \backslash j})$, the gradients on edge $i\to j$ are passed backward to its upstream
edges $\{k\to i\}_{k\in \mathcal{N}_{i} \backslash j}$ through the chain rule, as illustrated in the right panel of Fig.~\ref{fig:bilevel_flow_control_results}(a).
The full gradient on a non-targeted edge $k\to i$ can be obtained by summing the gradients on its downstream edges, computed as 
\begin{equation}
\frac{\partial\mathcal{O}}{\partial m_{k\to i}}=\sum_{l\in\mathcal{N}_{i}\backslash k}\,\,\sum_{m_{i\to l}\in\{\alpha_{i\to l},\hat{x}_{i\to l}\}}\frac{\partial\mathcal{O}}{\partial m_{i\to l}}\frac{\partial m_{i\to l}}{\partial m_{k\to i}}.
\end{equation}
The gradient messages $\{\frac{\partial\mathcal{O}}{\partial m_{k\to i}}\}$ are passed in a random and asynchronous manner, resulting in a decentralized algorithm. 

The gradient with respect to the control parameter on the non-targeted edge $(k,i)$ can be obtained straightforwardly as 
\begin{equation}
\frac{\partial\mathcal{O}}{\partial r_{ki}}=\sum_{m\in\{\alpha,\hat{x}\}}\bigg(\frac{\partial\mathcal{O}}{\partial m_{k\to i}}\frac{\partial m_{k\to i}}{\partial r_{ki}}+\frac{\partial\mathcal{O}}{\partial m_{i\to k}}\frac{\partial m_{i\to k}}{\partial r_{ki}}\bigg),
\end{equation}
which serves to update the control parameter in a gradient descent manner $r_{ki}\leftarrow r_{ki}-s\frac{\partial\mathcal{O}}{\partial r_{ki}}$ with certain step size $s$. 
The gradient for targeted edges can be similarly defined~\cite{BoLi2021sup}.
The control parameters are bounded to be $r_{ij}\in[0.9,1.1]$, achieved by necessary thresholding after gradient descent updates.
In this flow model, the gradient $\frac{\partial\mathcal{O}}{\partial r_{ki}}$ can be calculated exactly, leading to a global gradient descent (GGD) algorithm. 
However, the GGD approach requires computing the inverse of the Laplacian matrix in every iteration, which can be time-consuming for large networks. On the contrary, the gradients are computed in a local and distributed manner in the MP approach.
Similar ideas of gradient propagation of MP have been proposed in Refs~\cite{Eaton2009,Domke2013} in the context of approximate inference, which are usually implemented centrally in the reversed order of MP updates, unlike the decentralized approach presented here. 

The gradient computed by the MP algorithm provides an excellent estimation to the exact gradient, as illustrated in Fig.~\ref{fig:bilevel_flow_control_results}(b). 
For bilevel optimization, we do not wait for the convergence of the gradient-passing, but update the control parameters during the MP iterations to make the algorithm more efficient. It provides approximated gradient information, which is already effective for optimizing the global objective, as shown in Fig.~\ref{fig:bilevel_flow_control_results}(c).
The MP approach yields similar success rates in managing the network flows for different thresholds compared to the GGD approach as shown in Fig.~\ref{fig:bilevel_flow_control_results}(d), demonstrating the effectiveness of the MP approach for the bilevel optimization.

In summary, we propose MP algorithms for solving bilevel optimization in flow networks, focusing on applications in the routing game and flow control problems. In routing games, the objective functions in both levels admit a similar structure, which leads to two sets of similar messages being passed. 
Updates of the control variables based on localized information appear effective for toll optimization. 
However, the long-range impact of control-variable changes should be considered in some applications. 
This is accommodated by a separate distributed gradient-passing process, which is effective and efficient in flow control problems. 
Leveraging the sparse network structure, the MP approach offers efficient and intrinsically distributed algorithms in contrast to global optimization methods such as nonlinear programming, 
which is more generic, but is generally not scalable and therefore unsuitable for large-scale systems.
The MP approach  provides effective algorithms for bilevel optimization problems that are intractable or difficult to solve by global optimization approaches, such as combinatorial problems.
We believe that these MP methods provide a valuable tool for solving difficult bilevel optimization problems, especially in systems with sparsely coupled structures.

Source codes of this work can be found in \url{https://github.com/
boli8/bilevelMP_flow}.

\begin{acknowledgments}
We thank K. Y. Michael Wong, Tat Shing Choi, and Ho Fai Po for helpful discussions. B.L. and D.S. acknowledge support from the Leverhulme Trust (RPG-2018-092), European Union's Horizon 2020 research and innovation programme under the Marie Sk{\l}odowska-Curie Grant Agreement No. 835913.
B.L. acknowledges support from the startup funding from Harbin Institute of Technology, Shenzhen (Grant No. 20210134).
D.S. acknowledges support from the EPSRC programme grant TRANSNET (EP/R035342/1). 
C.H.Y. is supported by the Research Grants Council of the Hong Kong Special Administrative Region, China (Projects No. EdUHK GRF 18304316, No. GRF 18301217, and No. GRF 18301119), the Dean's Research Fund of the Faculty of Liberal Arts and Social Sciences (Projects No. FLASS/DRF 04418, No. FLASS/ROP 04396, and No. FLASS/DRF 04624), and the Internal Research Grant (Project No. RG67 2018-2019R R4015 and No. RG31 2020-2021R R4152), The Education University of Hong Kong, Hong Kong Special Administrative
Region, China.

\end{acknowledgments}

\bibliography{reference}

\end{document}


\title{Bilevel Optimization in Flow Networks -- A Message-passing Approach
\\ -- Supplemental Material}

\author{Bo Li}
\affiliation{Non-linearity and Complexity Research Group, Aston University, Birmingham,
B4 7ET, United Kingdom}
\affiliation{School of Science, Harbin Institute of Technology (Shenzhen), Shenzhen, 518055, China}

\author{David Saad}
\affiliation{Non-linearity and Complexity Research Group, Aston University, Birmingham,
B4 7ET, United Kingdom}

\author{Chi Ho Yeung}
\affiliation{Department of Science and Environmental Studies, The Education University of Hong Kong, 10 Lo Ping Road, Taipo, Hong Kong}

\maketitle

\section{Message-passing Algorithms For Non-atomic Routing Games\label{sec:MP_non_atomic}}

In this section, we provide details of the message-passing (MP) algorithm for non-atomic routing games in road networks, modeled by a directed graph $G(V,E)$. We denote $\mathbb{R}, \mathbb{Z}, \mathbb{N}$ as the sets of real numbers, integers and non-negative integers, respectively.

\subsection{Problem Setting and Notation for Non-atomic Games\label{subsec:def_nonatomic}}

A directed edge $e$ in the directed graph is represented by an ordered tuple $e=(i,j)$, where node $i$ is the head and node $j$ is the tail of edge $e$, i.e., $i=h(e),j=t(e)$. 
Note that there can be at most two directed edges connecting node $i$ and node $j$, i.e., $e=(i,j),e'=(j,i)$.

We write the set of incoming edges to node $i$ as $\partial_{i}^{\text{in}}=\{e|e\in E,t(e)=i\}$,
the set of outgoing edges from node $i$ as $\partial_{i}^{\text{out}}=\{e|e\in E,h(e)=i\}$
and the set of edges adjacent to node $i$ as $\partial i=\partial_{i}^{\text{in}}\cup\partial_{i}^{\text{out}}$.
For convenience, we define the incident operator $B:E\to V$,
with matrix elements
\begin{equation}
B_{i,e}=\begin{cases}
1, & \text{if }e\in\partial_{i}^{\text{in}}\\
-1, & \text{if }e\in\partial_{i}^{\text{out}}\\
0, & \text{otherwise.}
\end{cases}
\end{equation}

Consider the scenario where all users travel to a universal destination $\mathcal{D}$. The edge flows $\{ x_{e} \in \mathbb{R} \}$ resulting from users' path choices satisfy the flow conservation constraints,
\begin{align}
R_{i} & :=\Lambda_{i}+\sum_{e\in\partial_{i}^{\text{in}}}x_{e}-\sum_{e\in\partial_{i}^{\text{out}}}x_{e}\nonumber \\
 & =\Lambda_{i}+\sum_{e\in\partial i}B_{i,e}x_{e}=0,\quad\forall i\neq\mathcal{D},\label{eq:flow_conservation}
\end{align}
and the non-negativity constraints
\begin{equation}
x_{e}\geq0,\quad\forall e.\label{eq:flow_nonnegativity}
\end{equation}
Due to the flow conservation constraint, any resource on a leaf node
$i$ with only one outgoing edge (i.e., $|\partial_{i}^{\text{out}}|=1,|\partial_{i}^{\text{in}}|=0$)
must be transmitted to its only neighboring node $j$. Similarly,
if a leaf node $i$ with only one incoming edge (i.e., $|\partial_{i}^{\text{in}}|=1,|\partial_{i}^{\text{out}}|=0$)
is the destination node, then traffic must first arrive at its
only neighboring node $j$, and then go through the edge $(j,i)$ to the destination. In the former case, one can remove the leaf node
$i$ and add $\Lambda_{i}$ resources to its neighboring node $j$.
In the latter, one can simply set node $j$ as the destination.
By preprocessing the network using the above reduction, we can reduce the network to have no leaf nodes.

Denoting $\ell_{e}(x_{e})$ as the latency function on edge $e$ (assumed to be a non-decreasing function of $x_{e}$) and $\tau_{e}$ to be the corresponding toll,
the Wardrop equilibrium can be obtained by minimizing the following
potential function
\begin{equation}
\Phi(\boldsymbol{x})=\sum_{e\in E}\int_{0}^{x_{e}} \big[ \ell_{e}(y) + \tau_{e} \big] \mathrm{d}y=:\sum_{e\in E}\phi_{e}(x_{e}),\label{eq:Phi_def_nonatomic}
\end{equation}
subject to the flow conservation Eq.~(\ref{eq:flow_conservation}) and non-negativity constraints Eq.~(\ref{eq:flow_nonnegativity}).
We have assumed the same gauge between latency and toll can be used for all
users (more precisely the edge cost for a user is
$\ell_{e}(x_{e})+\chi\tau_{e}$
with $\chi$ being a coefficient converting money to time
which is set to one in Eq.~(\ref{eq:Phi_def_nonatomic}) in some appropriate unit).

The social cost of the routing game is defined as 
\begin{equation}
H(\boldsymbol{x})=\sum_{e\in E}x_{e}\ell_{e}(x_{e})=:\sum_{e\in E}\sigma_{e}(x_{e}),
\end{equation}
where the corresponding minimizer is the social optimum. Tolls are not assumed to contribute to the social cost $H(\boldsymbol{x})$.

The network planners only need to know the aggregated network flow on each edge and the corresponding latency in order to determine the social cost and set the tolls, while the specific paths where users choose are not that relevant. For this reason, we do not address the problem of finding individual routes for each user. The individual routing problem can be tackled by a multi-commodity formalism~\cite{Lonardi2021, Bonifaci2020}, or by using some physics-inspired algorithms~\cite{Yeung2013}.

\subsection{Intuition of the potential function $\Phi(\boldsymbol{x})$}

\begin{figure}
    \centering
    \includegraphics[scale=0.25]{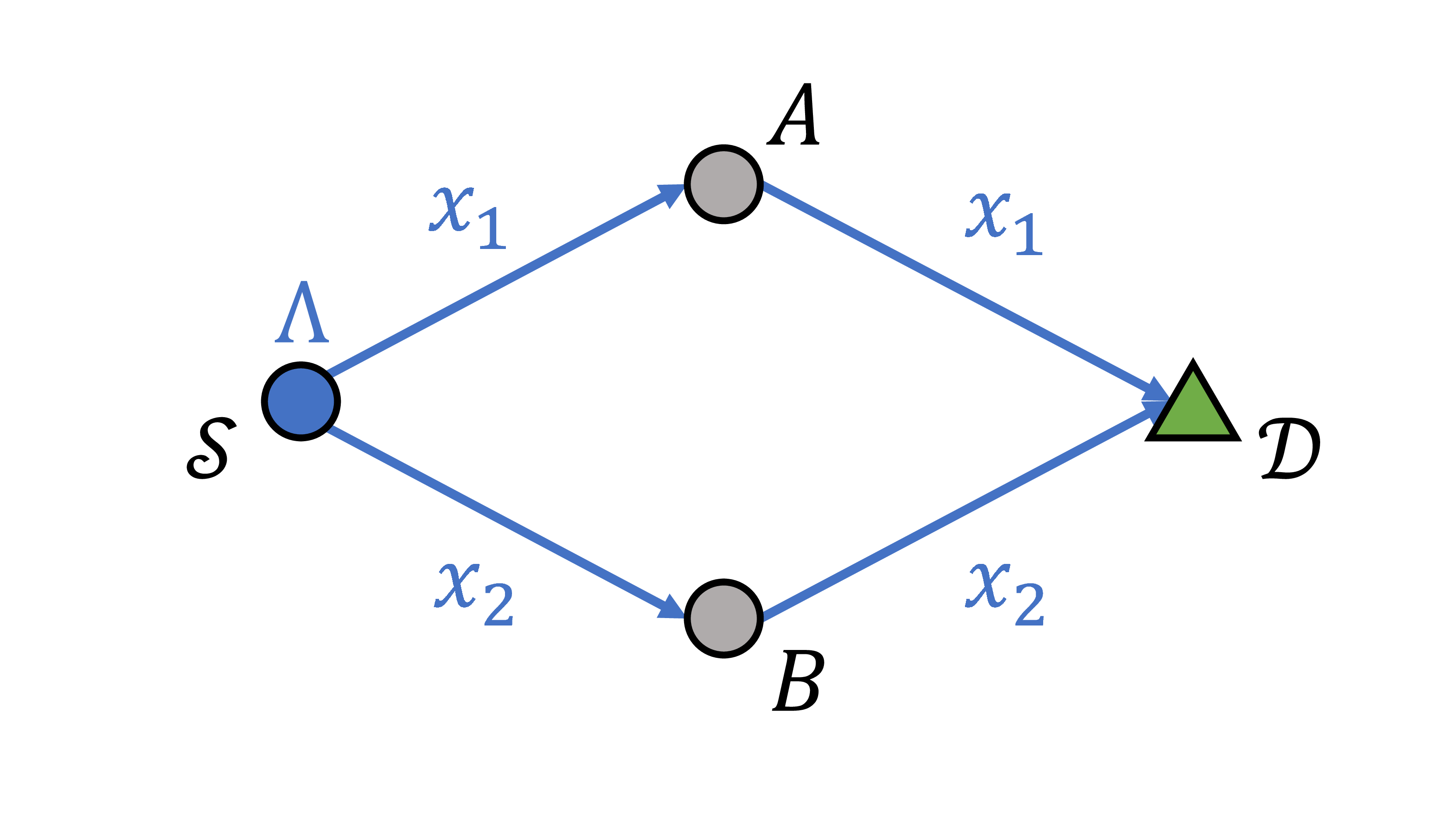}
    \caption{A small network with 4 nodes. There are $\Lambda> 0$ (where $\Lambda\in\mathbb{R}$ ) users originating from the source node $\mathcal{S}$ to the destination node $\mathcal{D}$, where $x_{1}$ users choose the path $\mathcal{P}_{1} = ((\mathcal{S}, A),(A, \mathcal{D}))$ and $x_{2}$ users choose the path $\mathcal{P}_{2} = ((\mathcal{S}, B),(B, \mathcal{D}))$.}
    \label{fig:4node}
\end{figure}

To gain some intuition of the role of the potential function $\Phi(\boldsymbol{x})$ in finding the Nash equilibrium, we consider the following simple scenario as shown in Fig.~\ref{fig:4node}, where there are $\Lambda > 0$ amount of users originating from the source node $\mathcal{S}$ to the destination node $\mathcal{D}$, and there are only two non-overlapping paths from node $\mathcal{S}$ to node $\mathcal{D}$, which are $\mathcal{P}_1 = ((\mathcal{S}, A),(A, \mathcal{D}))$ and $\mathcal{P}_{2} = ((\mathcal{S}, B),(B, \mathcal{D}))$. Denote $x_{1}$ and $x_{2}$ as the flow on path $\mathcal{P}_{1}$ and path $\mathcal{P}_2$ respectively. The flow conservation constraint asserts that $\Lambda = x_{1} + x_{2}$.

Users choosing $\mathcal{P}_{1}$ experience a cost $C_{1}(x_{1}) = \ell_{\mathcal{S} A}(x_{1}) + \tau_{\mathcal{S} A} + \ell_{A \mathcal{D}}(x_{1}) + \tau_{A \mathcal{D}}$, while users choosing $\mathcal{P}_{2}$ experience a cost $C_{2}(x_{2}) = \ell_{\mathcal{S} B}(x_{2}) + \tau_{\mathcal{S} B} + \ell_{B \mathcal{D}}(x_{2}) + \tau_{B \mathcal{D}}$. There are 3 possible scenarios of the Wardrop (Nash) equilibrium $x^{*}_{1}, x^{*}_{2}$, depending on the network parameters:
\begin{align}
\text{Case I:} & \qquad C_{1}(x^{*}_{1}) < C_{2}(x^{*}_{2}), \qquad x^{*}_{1} > 0, x^{*}_{2} = 0, \nonumber \\
\text{Case II:} & \qquad C_{1}(x^{*}_{1}) > C_{2}(x^{*}_{2}), \qquad x^{*}_{1} = 0, x^{*}_{2} > 0, \label{eq:Nash_eq_4node} \\
\text{Case III:} & \qquad C_{1}(x^{*}_{1}) = C_{2}(x^{*}_{2}), \qquad x^{*}_{1} > 0, x^{*}_{2} > 0. \nonumber
\end{align}
In Case I, all users choose $\mathcal{P}_{1}$ and there is no incentive for them to move to $\mathcal{P}_{2}$ as the corresponding cost $C_{2}(x_{2}=0)$ is higher. Similar analysis applies to Case II. In Case III, a user (controlling $\mathrm{d}x$ amount of traffic) choosing $\mathcal{P}_{1}$ also has no incentive to switch to $\mathcal{P}_{2}$; if she did so, her cost will become $C_{2}(x_{2} + \mathrm{d}x) \geq C_{2}(x_{2}) = C_{1}(x_{1})$, which is unfavorable.

Now we turn to the optimization problem as stated above
\begin{align}
\min_{x_{1}, x_{2}}  \Phi(x_{1}, x_{2}) & = \int_{0}^{x_{1}} \big[ \ell_{\mathcal{S}A}(y) + \tau_{\mathcal{S}A} \big] \mathrm{d}y + \int_{0}^{x_{1}} \big[ \ell_{A\mathcal{D}}(y) + \tau_{A\mathcal{D}} \big] \mathrm{d}y \\ 
& + \int_{0}^{x_{2}} \big[ \ell_{\mathcal{S}B}(y) + \tau_{\mathcal{S}B} \big] \mathrm{d}y + \int_{0}^{x_{2}} \big[ \ell_{B\mathcal{D}}(y) +\tau_{B\mathcal{D}} \big] \mathrm{d}y \nonumber \\
\text{s. t. } x_{1} & \geq 0, x_{x} \geq 0, \nonumber \\ 
& x_{1}+x_{2}=\Lambda, \nonumber
\end{align}
which can be solved by extremizing the Lagrangian function $\mathcal{L} = \Phi(x_{1}, x_{2}) - \mu (x_{1} + x_{2} - \Lambda) - \nu_{1} x_{1} - \nu_{2} x_{2}$, together with the Karush-Kuhn-Tucker (KKT) conditions $\nu_{1}, \nu_{2} \geq 0, \nu_{1} x_{1} = 0, \nu_{2} x_{2} = 0$. The extremum $\frac{\mathrm{d} \mathcal{L}}{\mathrm{d}x_{1}} = \ell_{\mathcal{S}A}(x_{1}) + \tau_{\mathcal{S}A} + \ell_{A\mathcal{D}}(x_{1}) + \tau_{A\mathcal{D}} - \mu - \nu_{1} = 0$ yields $C_{1}(x_{1}) = \mu + \nu_{1}$. Similarly, $C_{2}(x_{2}) = \mu + \nu_{2}$. The optimal dual parameters $(\mu^{*}, \nu^{*}_{1}, \nu_{2})$ has 3 possible scenarios under the KKT condition
\begin{align}
\text{Case I:} & \qquad \nu^{*}_{1} = 0, \nu^{*}_{2} > 0, \nonumber \\
\text{Case II:} & \qquad \nu^{*}_{1} > 0, \nu^{*}_{2} = 0, \label{eq:opt_Phi_4node} \\
\text{Case III:} & \qquad \nu^{*}_{1} = 0, \nu^{*}_{2} = 0, \nonumber
\end{align}
which exactly corresponds to the Wardrop equilibrium stated in Eq.~(\ref{eq:Nash_eq_4node}). Therefore, we have established that in this simple example, the Wardrop equilibrium can be identified by minimizing $\Phi(\boldsymbol{x})$.

We would also like to emphasize that the potential function $\Phi(\boldsymbol{x})$ does not carry any meaningful information about the routing costs (of either individuals of their aggregation), which can be seen by its definition. Therefore, the magnitude of $\Phi(\boldsymbol{x})$ is \emph{not} a performance measure of the traffic network. It plays an auxiliary role in defining the equilibrium flow as $\boldsymbol{x}^* = \text{arg}\min_{\boldsymbol{x}} \Phi(\boldsymbol{x}) \text{ s.t. }$ Eq.~(\ref{eq:flow_conservation}) and Eq.~(\ref{eq:flow_nonnegativity}). 
Roughly speaking, in the lower-level optimization problem, only $\text{arg}\min_{\boldsymbol{x}} \Phi(\boldsymbol{x})$ carries a correspondence to the routing problem, but not $\min_{\boldsymbol{x}} \Phi(\boldsymbol{x})$.

On the other hand, the social cost $H(\boldsymbol{x})$ corresponds to the aggregated travel latency of all users, which is a useful measure of the routing condition on a road network.

\subsection{MP Equations for Smooth Message Functions\label{subsec:MP_smooth}}

The MP equation for minimizing the potential $\Phi(\boldsymbol{x})$ reads
\begin{align}
\Phi_{i\to e}(x_{e}) & =\min_{\{x_{e'}\ge0\}|R_{i}=0}\sum_{e'\in\partial i\backslash e}\bigg[\Phi_{e'\to i}(x_{e'})+\phi_{e'}(x_{e'})\bigg],\label{eq:cavity_Phi_ie}
\end{align}
where the message function $\Phi_{i\to e}(x_{e})$ is called the cavity
energy in the jargon of statistical physics. Denoting $e=(i,j)$ and
$e'=(k,i)$, we can write $\Phi_{e'\to i}(x_{e'})=\Phi_{k\to e'}(x_{e'})$. 

In this framework, one needs to keep track of the profile of the message
functions $\Phi_{i\to e}(x_{e})$, which is only practical if they
are restricted to a certain family of functions and are easy to
optimize. One can approximate the message function $\Phi_{i\to e}(x_{e})$
by its series expansion around the working point $\tilde{x}_{i\to e}$~\cite{Wong2007}
\begin{align}
\Phi_{i\to e}(x_{e}) & =\Phi_{i\to e}(\tilde{x}_{i\to e}+\varepsilon_{e})\nonumber \\
 & \approx\text{\ensuremath{\Phi_{i\to e}}(\ensuremath{\tilde{x}_{i\to e}}})+\beta_{i\to e}|_{\tilde{x}_{i\to e}}\cdot\varepsilon_{e}+\frac{1}{2}\alpha_{i\to e}|_{\tilde{x}_{i\to e}}\cdot\big(\varepsilon_{e}\big)^{2},
\end{align}
where $\beta_{i\to e}$ and $\alpha_{i\to e}$ are the first and second
derivatives of $\Phi_{i\to e}$ evaluated at the working point $\tilde{x}_{i\to e}$,
assuming the message function $\Phi_{i\to e}(x_{e})$ is smooth in
the vicinity of $\tilde{x}_{i\to e}$. The MP equations have been
derived in~\cite{Wong2007} for undirected flow networks. Here,
we extend it to directed graph with non-negativity flow constraints.

Similarly, the interaction term $\phi_{e'}(x_{e'})$ is also approximated
as $\phi_{e'}(x_{e'})\text{\ensuremath{\approx\phi_{e'}}(\ensuremath{\tilde{x}_{i\to e}}})+\phi_{e'}'(\tilde{x}_{k\to e'})\varepsilon_{e}+\frac{1}{2}\phi_{e'}''(\tilde{x}_{k\to e'})\big(\varepsilon_{e}\big)^{2}$.
To solve the local optimization problem in Eq.~(\ref{eq:cavity_Phi_ie})
over the variables on edges $\{k\to e'|e'\in\partial i\backslash e\}$,
we introduce the Lagrangian
\begin{align}
L_{i\to e}= & \sum_{e'\in\partial i\backslash e}\bigg[\frac{1}{2}\alpha_{k\to e'}\big(\varepsilon_{e'}\big)^{2}+\beta_{k\to e'}\varepsilon_{e'}+\frac{1}{2}\phi_{e'}''(\tilde{x}_{k\to e'})\big(\varepsilon_{e'}\big)^{2}+\phi_{e'}'(\tilde{x}_{k\to e'})\varepsilon_{e'}\bigg]\nonumber \\
 & +\mu_{i\to e}R_{i}+\sum_{e'\in\partial i\backslash e}\lambda_{e'}(\tilde{x}_{k\to e'}+\varepsilon_{e'}),
\end{align}
where $\mu_{i\to e}$ and $\lambda_{e'}$ are the Lagrange multipliers
for the flow conservation constraint $R_{i}=0$ and flow non-negativity
constraint $x_{e'}\geq0$, respectively. Solving the extremum equation
$\frac{\partial L_{i\to e}}{\partial\varepsilon_{e'}}=0$ gives 
\begin{equation}
\varepsilon_{e'}^{*}(\mu_{i\to e})=\max\bigg(\frac{-1}{\alpha_{k\to e'}+\phi''_{e'}}\big(\mu_{i\to e}B_{i,e'}+\phi_{e'}'+\beta_{k\to e'}\big),-\tilde{x}_{k\to e'}\bigg),
\end{equation}
and the corresponding optimal cavity flow is
\begin{align}
x_{k\to e'}^{*}(\mu_{i\to e})=\tilde{x}_{k\to e'}+\varepsilon_{e'}^{*}(\mu_{i\to e}) & =\max\bigg(\tilde{x}_{k\to e'}-\frac{\mu_{i\to e}B_{i,e'}+\phi'_{e'}+\beta_{k\to e'}}{\alpha_{k\to e'}+\phi''_{e'}},0\bigg).\label{eq:x_kep_opt}
\end{align}
The Lagrange multiplier (or the dual variable) $\mu_{i\to e}$ needs to satisfy
\begin{equation}
R_{i\to e}(\mu_{i\to e};x_{e}):=\sum_{e'\in\partial i\backslash e}B_{i,e'}x_{k\to e'}^{*}(\mu_{i\to e})+B_{i,e}x_{e}+\Lambda_{i}=0.\label{eq:R_ie_of_mu}
\end{equation}

The function $R_{i\to e}(\mu;x_{e})$ is a non-increasing
piece-wise linear function of $\mu$. To determine the value of $\mu_{i\to e}^{*}$
at the optimum, we need to find the root of $R_{i\to e}(\mu;x_{e})$,
which can be done in finite steps by following the breakpoints of
the piece-wise linear function $R_{i\to e}(\mu;x_{e})$. Upon obtaining
the optimal dual variable $\mu_{i\to e}^{*}$, the messages $\beta_{i\to e}$
and $\alpha_{i\to e}$ are calculated by 
\begin{equation}
\beta_{i\to e}=\frac{\partial\Phi_{i\to e}^{*}(x_{e})}{\partial x_{e}}=\frac{\partial L_{i\to e}^{*}}{\partial x_{e}}=B_{i,e}\mu_{i\to e}^{*},\label{eq:beta_ie_MP}
\end{equation}
\begin{align}
\alpha_{i\to e} & =\frac{\partial^{2}\Phi_{i\to e}^{*}(x_{e})}{\partial x_{e}^{2}}=B_{i,e}\frac{\partial\mu_{i\to e}^{*}}{\partial x_{e}}=B_{i,e}\bigg[\frac{\partial x_{e}}{\partial\mu}\bigg|_{\mu=\mu_{i\to e}^{*}}\bigg]^{-1}\nonumber \\
 & =-\bigg[\frac{\partial}{\partial\mu}\sum_{e'\in\partial i\backslash e}B_{i,e'}x_{k\to e'}^{*}(\mu)\bigg|_{\mu=\mu_{i\to e}^{*}}\bigg]^{-1}\nonumber \\
 & =\bigg[\sum_{e'\in\partial i\backslash e}\frac{1}{\alpha_{k\to e'}+\phi''_{e'}}\Theta\bigg(\big(\alpha_{k\to e'}+\phi''_{e'}\big)\tilde{x}_{k\to e'}-\big(\mu_{i\to e}^{*}B_{i,e'}+\phi'_{e'}+\beta_{k\to e'}\big)\bigg)\bigg]^{-1},\label{eq:alpha_ie_MP}
\end{align}
where $\Theta(\cdot)$ is the Heaviside step function. The shadow
price interpretation of Lagrangian multiplier has been used in Eq.
(\ref{eq:beta_ie_MP}) and the inverse function theorem has been used
in Eq.~(\ref{eq:alpha_ie_MP}). In the implementation of the algorithm,
we take $x_{e}=\tilde{x}_{i\to e}$ in solving Eq.~(\ref{eq:R_ie_of_mu}).

\subsubsection{Destination node $\mathcal{D}$\label{subsec:treatmen_of_destination}}

There are two ways to treat the destination node $\mathcal{D}$:
\begin{itemize}
\item Method I: Since the destination node $\mathcal{D}$ has no constraint,
it will absorb all incoming flows (like a grounded
node in an electric circuit). So it has no preference for network flows of the incident edges, such that $\Phi_{\mathcal{D}\to e}(x_{e})=0$
and
\begin{equation}
\alpha_{\mathcal{D}\to e}=0,\quad\beta_{\mathcal{D}\to e}=0.\label{eq:alpha_beta_at_D}
\end{equation}
\item Method II: Alternatively, one can set an explicit constraint on the
flows to the destination node $\mathcal{D}$
\begin{equation}
R_{\mathcal{D}}:=\Lambda_{\mathcal{D}}+\sum_{e\in\partial\mathcal{D}}B_{\mathcal{D},e}x_{e}=0,
\end{equation}
where $\Lambda_{\mathcal{D}}=-\sum_{i\neq\mathcal{D}}\Lambda_{i}$.
Then the messages from the destination node $\mathcal{D}$ are calculated
in the same way as other nodes given by Eqs. (\ref{eq:beta_ie_MP})
and (\ref{eq:alpha_ie_MP}).
\end{itemize}
Method I is used for the experiments of routing games in the main
text.

\subsubsection{Working Points}

We also need a scheme to update the the working points $\{\tilde{x}_{i\to e}\}$
at which the messages $\{\alpha_{i\to e},\beta_{i\to e}\}$ are defined.
Here, we suggest to update the working point $\tilde{x}_{i\to e}$
such that it gets closer to the equilibrium flow $x_{e}^{*}$~\cite{Wong2007}
\begin{align}
x_{e}^{*} & =\text{arg}\min_{x_{e}\geq0}\big[\Phi_{i\to e}(x_{e})+\Phi_{j\to e}(x_{e})+\phi_{e}(x_{e})\big]\nonumber \\
 & =\text{arg}\min_{x_{e}\geq0}\bigg[\frac{1}{2}\big(\alpha_{i\to e}+\frac{1}{2}\phi''_{e}(\tilde{x}_{i\to e})\big)\big(x_{e}-\tilde{x}_{i\to e}\big)^{2}+\big(\beta_{i\to e}+\frac{1}{2}\phi_{e}'(\tilde{x}_{i\to e})\big)\big(x_{e}-\tilde{x}_{i\to e}\big)\nonumber \\
 & \qquad+\frac{1}{2}\big(\alpha_{j\to e}+\frac{1}{2}\phi''_{e}(\tilde{x}_{j\to e})\big)\big(x_{e}-\tilde{x}_{j\to e}\big)^{2}+\big(\beta_{j\to e}+\frac{1}{2}\phi_{e}'(\tilde{x}_{j\to e})\big)\big(x_{e}-\tilde{x}_{j\to e}\big)\bigg]\nonumber \\
 & =\max\bigg(\frac{\big(\alpha_{i\to e}+\frac{1}{2}\phi''_{i\to e}\big)\tilde{x}_{i\to e}+\big(\alpha_{j\to e}+\frac{1}{2}\phi''_{j\to e}\big)\tilde{x}_{j\to e}-(\beta_{i\to e}+\beta_{j\to e}+\frac{1}{2}\phi_{i\to e}'+\frac{1}{2}\phi_{j\to e}')}{\alpha_{i\to e}+\alpha_{j\to e}+\frac{1}{2}\phi''_{i\to e}+\frac{1}{2}\phi''_{j\to e}},0\bigg).
\end{align}
Furthermore, a learning rate $s$ is applied to update the working point
\begin{equation}
\tilde{x}_{i\to e}^{\text{new}}\leftarrow sx_{e}^{*}+(1-s)\tilde{x}_{i\to e}^{\text{old}},
\end{equation}
such that $\tilde{x}_{i\to e}$ does not jump too drastically; otherwise
the messages $\alpha_{i\to e}$ and $\beta_{i\to e}$ will approximate the curvature and slope of the message
function $\Phi_{i\to e}(x_{e})$ less precisely.

\subsection{Non-Smooth Message Functions}

\subsubsection{Qualitative Picture}

The MP algorithms in Sec.~\ref{subsec:MP_smooth} work well if the
smoothness assumption of the message function $\Phi_{i\to e}(x_{e})$
holds. However, it is not always the case in the routing game
problem, where the non-smoothness is induced by the non-negativity
constraints of Eq.~(\ref{eq:flow_nonnegativity}). Direct implementation
of the MP algorithms in Sec.~\ref{subsec:MP_smooth} leads to oscillations
of the messages when the traffic patterns are sparse. In fact, similar
non-convergence phenomena have been noticed in the system with a non-smooth
energy function~\cite{Wong2007}.

\begin{figure}
\includegraphics[scale=1.5]{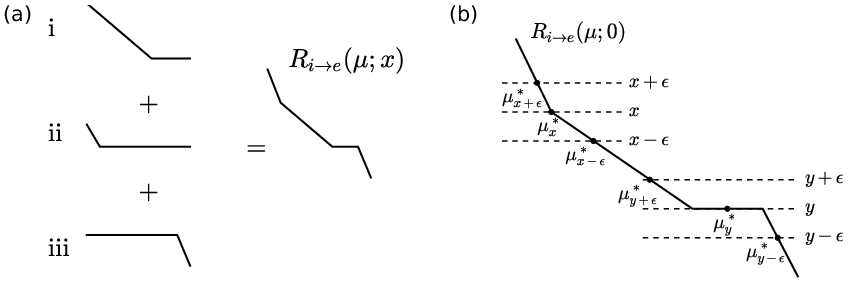} 

\caption{(a) The net resource $R_{i\to e}(\mu;x)$ defined in Eq.~(\ref{eq:R_ie_of_mu})
is a non-increasing piecewise-linear function of the Lagrange multiplier
$\mu$. Cases (i) and (ii) correspond to edges $e'$ incoming to node
$i$ ($B_{i,e'}=1$), while case (iii) corresponds to edge $e'$ outgoing
of node $i$ ($B_{i,e'}=-1$). (b) The roots of $R_{i\to e}(\mu;x_{e})$
in the vicinity of $x_{e}=x$ and $x_{e}=y$. It is assumed that edge
$e$ is an outgoing edge of node $i$ (with $B_{i,e}=-1$) such that
finding the root of $R_{i\to e}(\mu;x_{e})$ is equivalent to solving
$R_{i\to e}(\mu;0)=x_{e}$. For infinitesimal $\epsilon$, if the
flow $x_{e}$ changes from $x+\epsilon$ to $x-\epsilon$, the solution
of the Lagrange multiplier changes continuously from $\mu_{x+\epsilon}^{*}$
to $\mu_{x-\epsilon}^{*}$. On the other hand, there is a plateau
at $R_{i\to e}(\mu;0)=y$, so that when the flow $x_{e}$ changes
from $y+\epsilon$ to $y-\epsilon$, the solution of the Lagrange
multiplier changes discontinuously from $\mu_{y+\epsilon}^{*}$ to
$\mu_{y-\epsilon}^{*}$. \label{fig:R_vs_mu}}
\end{figure}

\begin{figure}
\includegraphics[scale=1.25]{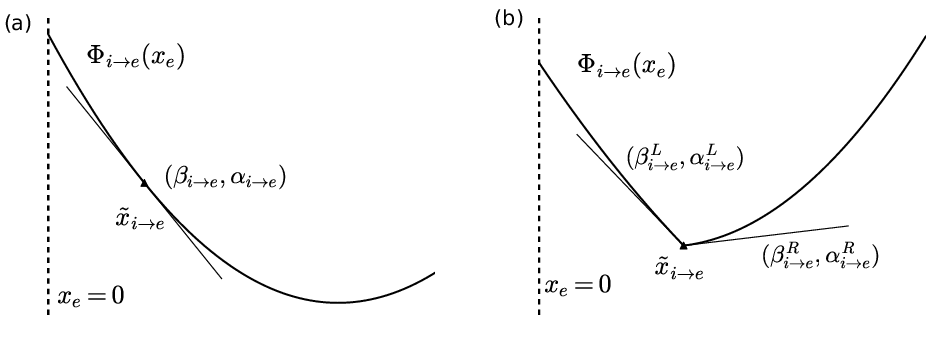}

\caption{(a) Smooth message function $\Phi_{i\to e}(x_{e})$, corresponding
to $x_{e}=x$ in Fig.~\ref{fig:R_vs_mu}(b). (b) Non-smooth message
function $\Phi_{i\to e}(x_{e})$ with one breakpoint, corresponding
to $x_{e}=y$ in Fig.~\ref{fig:R_vs_mu}(b), where the first and second
derivatives of $\Phi_{i\to e}(x_{e})$ are discontinuous near $x_{e}=y$.
\label{fig:cavity_energy_illustrate}}
\end{figure}

To better understand this phenomenon, we examine $R_{i\to e}(\mu;x_{e})$
as a function of the Lagrange multiplier $\mu$ in Eq.~(\ref{eq:R_ie_of_mu}),
of which the root $\mu^{*}$ (satisfying $R_{i\to e}(\mu^{*},x_{e})=0$)
will determine $\beta_{i\to e}$ and $\alpha_{i\to e}$ in Eqs.~(\ref{eq:beta_ie_MP})
and (\ref{eq:alpha_ie_MP}). The function $R_{i\to e}(\mu;x_{e})$
is non-increasing piecewise-linear function as illustrated in Fig.~\ref{fig:R_vs_mu}(a). Assuming edge $e$ is an outgoing edge of node
$i$ (with $B_{i,e}=-1$), finding the root of $R_{i\to e}(\mu;x_{e})$
is equivalent to solving $R_{i\to e}(\mu;0)=x_{e}$. Consider the
configuration in Fig.~\ref{fig:R_vs_mu}(b), where the solution of
$R_{i\to e}(\mu;0)=y$ occurs at a plateau, such that the solution
(denoted as $\mu_{y}^{*}$) is degenerate; when the flow $x_{e}$
changes infinitesimally from $y+\epsilon$ to $y-\epsilon$, the solution
of the Lagrange multiplier changes discontinuously from $\mu_{y+\epsilon}^{*}$
to $\mu_{y-\epsilon}^{*}$. In this case, the slope $\beta_{i\to e}$
of the cavity energy $\Phi_{i\to e}(x_{e})$ changes discontinuously
from $x_{e}=y+\epsilon$ to $x_{e}=y-\epsilon$, while the curvature
$\alpha_{i\to e}$ is ill-defined at $x_{e}=y$. The profiles of the
message function $\Phi_{i\to e}(x_{e})$ in the smooth and non-smooth
cases are illustrated in Fig.~\ref{fig:cavity_energy_illustrate}.
If the normal messages $\{\beta_{i\to e},\alpha_{i\to e}\}$ are used
when the message function $\Phi_{i\to e}(x_{e})$ is non-smooth,
the solution will be jumping between the two branches, resulting in
non-convergence behaviors of the MP algorithms as observed in Ref.~\cite{Wong2007}.

\subsubsection{Criteria for Non-smooth Message Function}

As mentioned above, the message function $\Phi_{i\to e}(x_{e})$ is non-smooth if the solution of $\mu$ in Eq.~(\ref{eq:R_ie_of_mu})
is degenerate. This occurs if the optimal flow $x_{k\to e'}^{*}(\mu_{i\to e})$
of all descendant edges $e'\in\partial i\backslash e$ are inactive,
i.e., lying in the zero branch of the function in Eq.~(\ref{eq:x_kep_opt});
when $B_{i,e}x_{e}+\Lambda_{i}=0$, the flow conservation equation
$R_{i\to e}(\mu;x_{e})=\sum_{e'\in\partial i\backslash e}B_{i,e'}x_{k\to e'}^{*}(\mu)=0$
has degenerate solutions. In this case, all the resources $\Lambda_{i}$
are transmitted along edge $e$, while the flows on all other edges
$\partial i\backslash e$ adjacent to node $i$ are idle. When $\Lambda_{i}>0$,
edge $i\to e$ is a leaf in the subgraph with edges holding non-zero
flows, therefore we call edge $i\to e$ a primary effective leaf in
such cases. Since $\Lambda_{i}\geq0$ and $x_{e}\geq0$, only the out-going edge
$i\to e$ from node $i$ (with $B_{i,e}=-1$) can be a primary
effective leaf.

\begin{figure}
\includegraphics[scale=0.8]{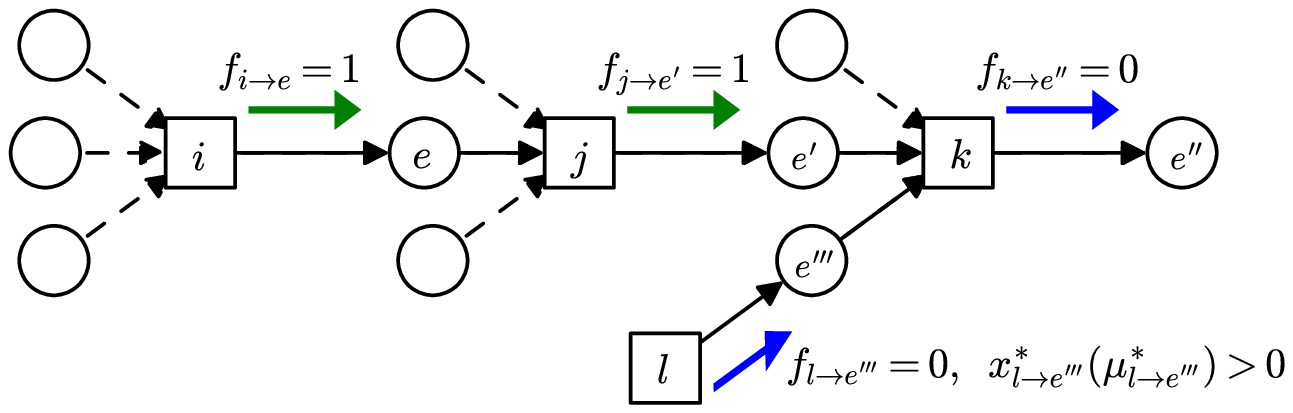}

\caption{Illustration of the effective leaf edges. Arrows with dashed lines
correspond to edges with zero optimal cavity flow $x^{*}(\mu^{*})$
in the MP calculation (expression given in Eq.~(\ref{eq:x_kep_opt}))
of one of its downstream edges. Edge $i\to e$ is a primary effective
leaf (assuming $\Lambda_{i}>0$), as all the upstream optimal flows
are zero. Edge $j\to e'$ is a general effective leaf, as its upstream
edges are either effective leaves or attain zero optimal cavity flows.
Edge $k\to e''$ is a non-effective leaf, as its upstream edge $l\to e'''$
is a non-effective leaf and it has a non-zero optimal flow $x_{l\to e'''}^{*}(\mu_{l\to e'''}^{*})>0$.
\label{fig:effective_leaf_illustrate}}
\end{figure}

The leaf state can also propagate from primary effective leaves to
downstream edges. We define edge $i\to e$ to be a general effective
leaf if and only if $\forall e'\in\partial i\backslash e$, either
(i) the optimal flow $x_{k\to e'}^{*}(\mu_{i\to e}^{*})=0$ in Eq.
(\ref{eq:x_kep_opt}) or (ii) edge $k\to e'$ is a general effective
leaf. A primary leaf is by default a general effective leaf. If an
edge $i\to e$ is an effective leaf, we denote $f_{i\to e}=1$, otherwise
$f_{i\to e}=0$. An example of effective leaf configurations is shown
in Fig \ref{fig:effective_leaf_illustrate}. 

It can be proved by contradiction that only out-going edges $i\to e$
with $B_{i,e}=-1$ can be general effective leaves under the condition
$\Lambda_{i}\geq0$. The set of effective leaves in the upstream of
edge $i\to e$ is $EL_{i\to e}=\{e'|e'\in\partial i\backslash e,f_{k\to e'}=1\}$,
while the set of non-effective leaves is $NEL_{i\to e}=\{e'|e'\in\partial i\backslash e,f_{k\to e'}=0\}$. 

Since there is at most one plateau in the function $R_{i\to e}(\mu;x)$,
the cavity message $\Phi_{i\to e}(x_{e})$ has at most one breakpoint.
For an effective leaf edge $i\to e$, we always use the breakpoint
of the message function (denoted as $\tilde{x}_{i\to e}^{b}$) as
the working point, such that $\Phi_{i\to e}(x_{e})$ has the following
expression
\begin{equation}
\Phi_{i\to e}(x_{e})=\begin{cases}
\frac{1}{2}\alpha_{i\to e}^{L}(x_{e}-\tilde{x}_{i\to e}^{b})^{2}+\beta_{i\to e}^{L}(x_{e}-\tilde{x}_{i\to e}^{b})+E_{i\to e}(\tilde{x}_{i\to e}^{b}) & x<\tilde{x}_{i\to e}^{b},\\
\frac{1}{2}\alpha_{i\to e}^{R}(x_{e}-\tilde{x}_{i\to e}^{b})^{2}+\beta_{i\to e}^{R}(x_{e}-\tilde{x}_{i\to e}^{b})+E_{i\to e}(\tilde{x}_{i\to e}^{b}) & x>\tilde{x}_{i\to e}^{b}.
\end{cases}\label{eq:non_smooth_cavity_Phi_ie}
\end{equation}
For a primary effective leaf edge $i\to e$, the breakpoint is $\tilde{x}_{i\to e}^{b}=\Lambda_{i}$.
For a general effective leaf edge $i\to e$, the breakpoint is most
likely (but not always) located at the value of effective resource
defined as 
\begin{equation}
\Lambda_{i\to e}^{\text{eff}}:=\Lambda_{i}+\sum_{e'\in EL_{i\to e}}B_{i,e'}\tilde{x}_{i\to e}^{b}.\label{eq:effective_Lambda}
\end{equation}

\subsection{MP Equations for Non-smooth Message Functions}\label{sec:MP_for_nonsmooth}

The MP equations for non-smooth message functions can be obtained
with the information on effective leaf status of upstream edges,
where one replaces the quadratic expansion $\Phi_{i\to e}(\tilde{x}_{i\to e}+\varepsilon_{e})\approx\text{\ensuremath{\Phi_{i\to e}}(\ensuremath{\tilde{x}_{i\to e}}})+\beta_{i\to e}\varepsilon_{e}+\frac{1}{2}\alpha_{i\to e}(\varepsilon_{e})^{2}$
by the piecewise quadratic counterpart in Eq.~(\ref{eq:non_smooth_cavity_Phi_ie})
when edge $i\to e$ is determined to be an effective leaf, and the double-sided
message parameters $\{\alpha_{i\to e}^{L},\beta_{i\to e}^{L},\alpha_{i\to e}^{R},\beta_{i\to e}^{R}\}$
are maintained and passed to its downstream edges.

For updating the messages, the computation of $\min_{\{x_{e'}\ge0\}|R_{i}=0}\sum_{e'\in\partial i\backslash e}\big[\Phi_{k\to e'}(x_{e'})+\phi_{e'}(x_{e'})\big]$
can be tedious if there are multiple effective leaf edges in $\{k\to e'|e'\in\partial i\backslash e\}$,
where one needs to solve for a quadratic optimization of every case,
where one branch of each non-smooth message function is selected each time (there
are $2^{|EL_{i\to e}|}$ such cases in total). To simplify this
process, we propose to firstly fix the flow $x_{e'}$ of effective
leaves $EL_{i\to e}$ to be their breakpoints $\tilde{x}_{k\to e'}^{b}$
and then optimize non-effective leaf edges $NEL_{i\to e}$
\begin{equation}
\min_{\{\varepsilon_{e'}|e'\in NEL_{i\to e}\}}\sum_{e'\in NEL_{i\to e}}\bigg[\Phi_{k\to e'}(\tilde{x}_{k\to e'}+\varepsilon_{e'})+\phi(\tilde{x}_{k\to e'}+\varepsilon_{e'})\bigg],\label{eq:min_Phi_NEL}
\end{equation}
\begin{align}
\text{s. t. }0 & =\sum_{e'\in NEL_{i\to e}}B_{i,e'}\big[\tilde{x}_{k\to e'}+\varepsilon_{e'}\big]+\sum_{e'\in EL_{i\to e}}B_{i,e'}\tilde{x}_{k\to e'}^{b}+B_{i,e}x_{e}+\Lambda_{i}\nonumber \\
 & =\sum_{e'\in NEL_{i\to e}}B_{i,e'}\big[\tilde{x}_{k\to e'}+\varepsilon_{e'}\big]+\Lambda_{i\to e}^{\text{eff}}+B_{i,e}x_{e},\\
0 & \leq\tilde{x}_{k\to e'}+\varepsilon_{e'}.
\end{align}

We then perturb the optimal solution by perturbing some of the flows
$x_{e'}$ of upstream edges $e'\in\partial i\backslash e$ by
an infinitesimal amount $\mathrm{d}x$ as $x_{e'}=x_{k\to e'}^{*}+\eta_{e'}\mathrm{d}x$
for non-effective leaves and $x_{e'}=\tilde{x}_{k\to e'}^{b}+\eta_{e'}\mathrm{d}x$
for effective leaves with $\eta_{e'}=0,\pm1$. For non-smooth message
function, $\eta_{e'}=-1$ and $\eta_{e'}=1$, corresponding to the
left and the right branch of $\Phi_{k\to e'}(x_{e'})$, respectively.
To obey the flow conservation constraint $R_{i}=0$, the perturbation
coefficient $\eta_{e'}$ must satisfy $\sum_{e'\in\partial i\backslash e}B_{i,e'}\eta_{e'}=0$.

If the perturbation configuration $\{\eta_{e'}^{*}\}$ leading to
the lowest energy of $\sum_{e'\in\partial i\backslash e}\big[\Phi_{e'\to i}(x_{e'})+\phi_{e'}(x_{e'})\big]$
reduces the outcome of Eq.~(\ref{eq:min_Phi_NEL}),
we need to consider adding the effective leafs $k\to e'$ with $\eta_{e'}^{*}\neq0$
as  active optimization variables in addition to the non-effective
leafs. Specifically, we define $\eta^{\text{active}}=\{e'|e'\in EL_{i\to e},\eta_{e'}^{*}\neq0\}$,
and proceed to solve 
\begin{equation}
\min_{\{\varepsilon_{e'}|e'\in NEL_{i\to e}\cup\eta^{\text{active}}\}}\sum_{e'\in NEL_{i\to e}\cup\eta^{\text{active}}}\bigg[\Phi_{k\to e'}(\tilde{x}_{k\to e'}+\varepsilon_{e'})+\phi(\tilde{x}_{k\to e'}+\varepsilon_{e'})\bigg],
\end{equation}
\begin{align}
\text{s. t. }0= & \sum_{e'\in NEL_{i\to e}\cup\eta^{\text{active}}}B_{i,e'}\big[\tilde{x}_{k\to e'}+\varepsilon_{e'}\big]\nonumber \\
 & \qquad+\sum_{e'\in EL_{i\to e}\backslash\eta^{\text{active}}}B_{i,e'}\tilde{x}_{k\to e'}^{b}+B_{i,e}x_{e}+\Lambda_{i},\\
0\leq & \tilde{x}_{k\to e'}+\varepsilon_{e'},
\end{align}
where we use $\Phi_{k\to e'}=\frac{1}{2}\alpha_{k\to e'}^{L}\varepsilon_{e'}^{2}+\beta_{k\to e'}^{L}\varepsilon_{e'}$
if $\eta_{e'}^{*}=-1$ and use $\Phi_{k\to e'}=\frac{1}{2}\alpha_{k\to e'}^{R}\varepsilon_{e'}^{2}+\beta_{k\to e'}^{R}\varepsilon_{e'}$
if $\eta_{e'}^{*}=1$.

The primal and dual variables in the optimum satisfy
\begin{align}
x_{k\to e'}^{*}(\mu)=\tilde{x}_{k\to e'}+\varepsilon_{e'}^{*}(\mu) & =\max\bigg(\tilde{x}_{k\to e'}-\frac{\mu B_{i,e'}+\phi'_{e'}+\beta_{k\to e'}}{\alpha_{k\to e'}+\phi''_{e'}},0\bigg),
\end{align}
\begin{align}
R_{i\to e}(\mu;x_{e})= & \sum_{e'\in NEL_{i\to e}\cup\eta^{\text{active}}}B_{i,e'}x_{k\to e'}^{*}(\mu)\nonumber \\
 & \qquad+\sum_{e'\in EL_{i\to e}\backslash\eta^{\text{active}}}B_{i,e'}\tilde{x}_{k\to e'}^{b}+B_{i,e}x_{e}+\Lambda_{i}=0.\label{eq:R_ie_of_mu_leaf}
\end{align}
If the solution $\mu^{*}$ in Eq.~(\ref{eq:R_ie_of_mu_leaf}) is non-degenerate, we have 
\begin{equation}
\beta_{i\to e}=B_{i,e}\mu^{*},\label{eq:beta_ie_MP_with_leaf}
\end{equation}
\begin{align}
\alpha_{i\to e}= & \bigg[\sum_{e'\in NEL_{i\to e}\cup\eta^{\text{active}}}\frac{1}{\alpha_{k\to e'}+\phi''_{e'}}\nonumber \\
 & \quad\times\Theta\bigg(\big(\alpha_{k\to e'}+\phi''_{e'}\big)\tilde{x}_{k\to e'}-\big(\mu^{*}B_{i,e'}+\phi'_{e'}+\beta_{k\to e'}\big)\bigg)\bigg]^{-1},\label{eq:alpha_ie_MP_with_leaf}
\end{align}
in which case the edge $i \to e$ is not an effective leaf with $f_{i \to e} = 0$.

On the other hand, if the solution $\mu^{*}$ is degenerate, we need
to consider $x_{e}=\tilde{x}_{i\to e}-\mathrm{d}x$ to solve for
$\beta_{i\to e}^{L},\alpha_{i\to e}^{L}$, and consider $x_{e}=\tilde{x}_{i\to e}+\mathrm{d}x$
to solve for $\beta_{i\to e}^{R},\alpha_{i\to e}^{R}$, and identify edge $i \to e$ as an effective leaf with $f_{i \to e} = 1$.

It can also be shown that $\beta_{i\to e}^{R}>\beta_{i\to e}^{L}$ and the non-smooth message function $\Phi_{i\to e}(x_{e})$ is convex.

\subsubsection{Update of the Working Points}

If the message function $\Phi_{i\to e}(x_{e})$ is non-smooth, we
would like to bring the working point $\tilde{x}_{i\to e}$ to the
vicinity of the breakpoint of the two branches. To determine whether
an edge $i\to e$ is an effective leaf, we perform the following procedure: 
We check the two following criteria: (i) each edge $k\to e'$ in the upstream
edge set $\partial i\backslash e$ satisfies either $f_{k\to e'}=1$
or $\tilde{x}_{k\to e'}=0$; (ii) the difference between the current
working point and the effective resource $|\tilde{x}_{i\to e}-\Lambda_{i\to e}^{\text{eff}}|$
is smaller than some threshold ($\Lambda_{i\to e}^{\text{eff}}$ is
defined in Eq.~(\ref{eq:effective_Lambda})). If both criteria (i)
and (ii) are met, then we use the effective resource as the working
point $\tilde{x}_{i\to e}=\Lambda_{i\to e}^{\text{eff}}$, and
perform the optimization $\min_{\{x_{e'}\ge0\}|R_{i}=0}\sum_{e'\in\partial i\backslash e}\big[\Phi_{k\to e'}(x_{e'})+\phi_{e'}(x_{e'})\big]$;
if it results in degenerate solutions of the Lagrangian multiplier
$\mu^{*}$ for the flow conservation constraint, then edge $i\to e$
is determined as an effective leaf and the double-sided messages $\{\beta_{i\to e}^{L},\alpha_{i\to e}^{L},\beta_{i\to e}^{R},\alpha_{i\to e}^{R}\}$
are computed. Otherwise, edge $i\to e$ is a non-effective leaf and
the normal messages $\{\beta_{i\to e},\alpha_{i\to e}\}$ are recorded. 

If criteria (i) and (ii) are not met, we use the current value of
the working point $\tilde{x}_{i\to e}$ to solve for the messages.
Similarly, if the optimization leads to degenerate solutions of $\mu^{*}$,
then edge $i\to e$ is determined as an effective leaf. Otherwise,
edge $i\to e$ is a non-effective leaf. 

Similar to the case of smooth message functions in Sec.~\ref{subsec:MP_smooth},
the working point is updated as
\begin{equation}
\tilde{x}_{i\to e}^{\text{new}}\leftarrow sx_{e}^{*}+(1-s)\tilde{x}_{i\to e}^{\text{old}},
\end{equation}
where $x_{e}^{*}=\text{arg}\min_{x_{e}\geq0}\big[\Phi_{i\to e}(x_{e})+\Phi_{j\to e}(x_{e})+\phi_{e}(x_{e})\big]$.

\subsection{Results of MP Algorithm for Routing Game}

The MP algorithm for solving the (single-level) equilibrium flow problem is summarized in Algorithm.~\ref{alg:MP_routing_lower}.

\begin{algorithm}
\SetNoFillComment
\caption{Message-passing algorithm for equilibrium flows in routing games (single-level, single destination, Method I to treat the destination node $\mathcal{D}$ )}\label{alg:MP_routing_lower}

\KwIn{Road network $G(V,E)$ (pre-processed to remove dangling nodes), node parameters $\{ \Lambda_i \}$ defining the resources, edge parameters defining the latency function $\ell_e(x_e)$ and the edge-wise potential $\phi_e(x_e)$ (defined in Eq.~(\ref{eq:Phi_def_nonatomic})), maximal number of iterations $T$.}

\BlankLine
Initialize the messages $\{ \alpha^{(m)}_{i \to e}, \beta^{(m)}_{i \to e} \}$ and $\{ f_{i \to e} \}$ (effective leaf states) randomly.

\BlankLine
Initialize the working points $\{ \tilde{x}_{i \to e} \}$ randomly.

\BlankLine
\For{$t \textnormal{ in } 1:T$}{

Randomly select a node $i$ and one of its adjacent edge $e \in \partial i$.

\BlankLine
\BlankLine
\textbf{Begin Subroutine (a) (update the messages $\alpha^{(m)}_{i \to e}, \beta^{(m)}_{i \to e}, f_{i \to e}$ using $\phi_e(\cdot)$):}

Compute $\Lambda^{\text{eff}}_{i \to e}$ using Eq.~(\ref{eq:effective_Lambda}).

\uIf{$i = \mathcal{D}$}{
Set $\alpha_{i \to e} = 0, \beta_{i \to e} = 0, f_{i \to e} = 0$.
}

\uElseIf{$\exists \; e' \in \{ \partial i \backslash e \}, \tilde{x}_{k \to e'} \neq 0, f_{k \to e'} = 0$ \textnormal{and} $|\tilde{x}_{i \to e} -  \Lambda^{\text{eff}}_{i \to e}| >$ \textnormal{threshold} }
{Update the messages $\alpha^{(m)}_{i \to e}, \beta^{(m)}_{i \to e}, f_{i \to e}$ with the working point evaluated at the up-to-date value of $\tilde{x}_{i \to e}$, following the procedures outlined in Sec.~\ref{sec:MP_for_nonsmooth}.
}

\Else{
\tcc{edge $i \to e$ is a potential effective leaf.}
Replace $\tilde{x}_{i \to e}$ by $\Lambda^{\text{eff}}_{i \to e}$ in Sec.~\ref{sec:MP_for_nonsmooth} to compute the messages $\alpha^{(m)}_{i \to e}, \beta^{(m)}_{i \to e}, f_{i \to e}$.

\If{$f_{i \to e} = 1$}{
\tcc{This is an optional step, corresponding to a more conservative strategy of identifying the effective leaf state.}
{Using the newly computed $\alpha^{(m)}_{i \to e}, \beta^{(m)}_{i \to e}$, check whether $\Lambda^{\text{eff}}_{i \to e}$ is the minimum of $\Phi^{\text{full}}(x_e)$. 

If yes, adopt the messages $\alpha^{(m)}_{i \to e}, \beta^{(m)}_{i \to e}, f_{i \to e}$.

If no, reset $f_{i \to e} = 0$}.
}

\If{$f_{i \to e} = 0$}{
\tcc{Do not consider edge $i \to e$ as an effective leaf.}
Use the up-to-date value of $\tilde{x}_{i \to e}$ to recompute the messages $\alpha^{(m)}_{i \to e}, \beta^{(m)}_{i \to e}$.
}

} 

\textbf{End Subroutine (a)}

\BlankLine
\BlankLine
\If{$f_{i \to e} = 1$}{
Set $\tilde{x}_{i \to e} = \Lambda^{\text{eff}}_{i \to e}$.
}

\BlankLine
\BlankLine

Using the up-to-date messages $\{ \alpha^{(m)}_{i \to e}, \beta^{(m)}_{i \to e}, f_{i \to e} \}$ and $\{ \tilde{x}_{i \to e} \}$ to determine the approximated forms of $\Phi_{i \to e}(x_e)$ and $\Phi_{j \to e}(x_e)$.

Compute $x_e^*=\arg \min _{x_e \geq 0}\left[\Phi_{i \to e}\left(x_e\right)+\Phi_{j \to e}\left(x_e\right)+\phi_e\left(x_e\right)\right]$.

Update the working point as $\tilde{x}_{i \to e}^{\text {new }} \leftarrow s x_e^*+(1-s) \tilde{x}_{i \to e}^{\text {old }}$ with a learning rate $s$.

\BlankLine
\BlankLine
\If{\textnormal{messages converge}}{
Exit the for loop.
}

} 

\BlankLine

For each edge $e$, compute the equilibrium flow as $x_e^*=\arg \min _{x_e \geq 0}\left[\Phi_{i \to e}\left(x_e\right)+\Phi_{j \to e}\left(x_e\right)+\phi_e\left(x_e\right)\right]$.

\BlankLine
\KwOut{Convergence status and equilibrium flows $\{ x^{*}_{e} \}$}
\end{algorithm}

\clearpage

\begin{figure}
\includegraphics[scale=1.5]{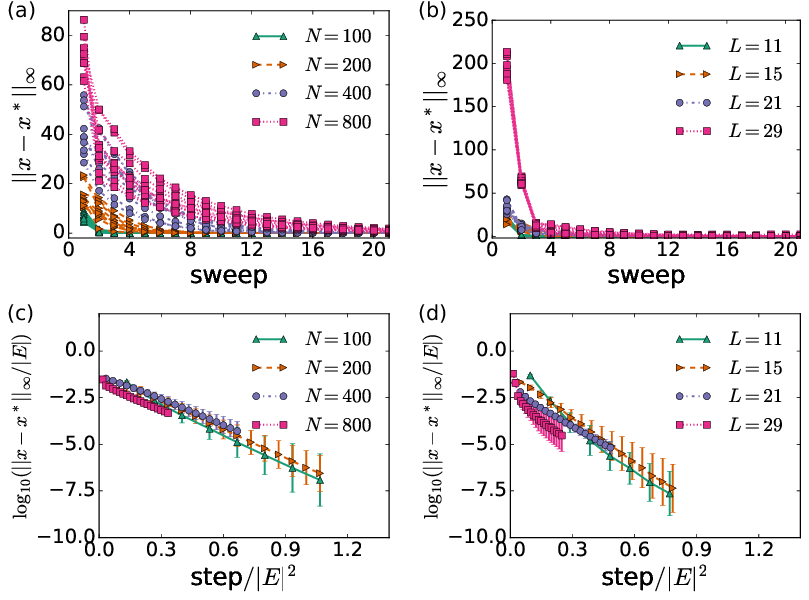}

\caption{Convergence of the MP algorithm for routing games in networks to the equilibrium flows $\boldsymbol{x}^{*}$. Random regular graphs
of degree $3$ are considered in (a)(c), while small-world networks
obtained by rewiring square lattices with randomly chosen shortcut
edges (rewiring probability $p_{rw}=0.05$) are considered in (b)(d), respectively.
The flows adjacent to the destination node $\mathcal{D}$ are unconstrained
(Method I in Sec.~\ref{subsec:treatmen_of_destination}), reminiscent
of a grounded node in electric circuits. In (a)(b), the MP process
for specific problem realizations is illustrated. Each sweep comprises $40|E|$
local MP updates. Panels (c)(d) are obtained by averaging
over $10$ problem realizations, and rescaling the MP runtime by $|E|^{2}$.
\label{fig:routing_game_Lower_MP_converge_with_grounded}}

\end{figure}

\begin{figure}
\includegraphics[scale=1.5]{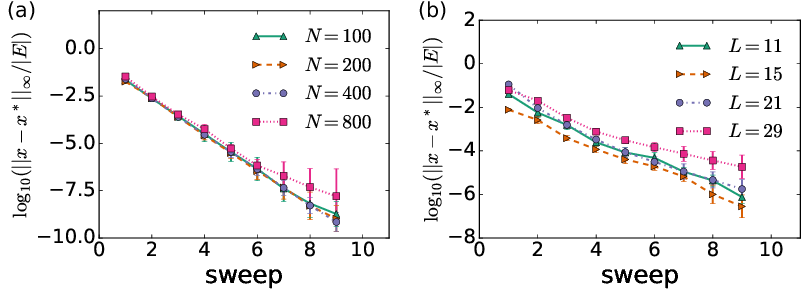}

\caption{Same setting as in Fig.~\ref{fig:routing_game_Lower_MP_converge_with_grounded},
except that the flows adjacent to the destination node $\mathcal{D}$
explicitly obey the constraint $R_{\mathcal{D}}=\Lambda_{\mathcal{D}}+\sum_{e\in\partial\mathcal{D}}B_{\mathcal{D},e}x_{e}=0$,
where $\Lambda_{\mathcal{D}}=-\sum_{i\protect\neq\mathcal{D}}\Lambda_{i}$
(Method II in Sec.~\ref{subsec:treatmen_of_destination}). Each sweep
comprises of $40|E|$ local MP updates. \label{fig:routing_game_Lower_MP_converge_no_grounded}}
\end{figure}

We further report results of the MP algorithms described
above. Taking into account the possible non-smooth structure of message
functions $\Phi_{i\to e}(x_{e})$, the MP algorithm converges well
for various types of graphs and resource distributions. We demonstrate
the effectiveness of the algorithm in Figures \ref{fig:routing_game_Lower_MP_converge_with_grounded}
and \ref{fig:routing_game_Lower_MP_converge_no_grounded}, where random
regular graphs and small-world networks are considered. The small-world
networks are obtained by rewiring square lattices with randomly chosen
shortcut edges~\cite{BoLiPRR2020}. In Fig.~\ref{fig:routing_game_Lower_MP_converge_with_grounded},
the flows adjacent to the destination node $\mathcal{D}$ are unconstrained
(Method I in Sec.~\ref{subsec:treatmen_of_destination}). The MP algorithms
converge to the correct equilibrium flows $\boldsymbol{x}^{*}$, and
the empirical complexity for computing the equilibrium flows up to
a certain error $|\boldsymbol{x}^{MP}-\boldsymbol{x}^{*}|$ is roughly
$O(|E|^{2})$.

In Fig.~\ref{fig:routing_game_Lower_MP_converge_no_grounded}, we
use Method II in Sec.\ref{subsec:treatmen_of_destination}, i.e.,
we put an explicit constraint to the flows adjacent to the destination
node $\mathcal{D}$ as 
\begin{equation}
R_{\mathcal{D}}=\Lambda_{\mathcal{D}}+\sum_{e\in\partial\mathcal{D}}B_{\mathcal{D},e}x_{e}=0,
\end{equation}
where $\Lambda_{\mathcal{D}}=-\sum_{i\neq\mathcal{D}}\Lambda_{i}$.
In this approach, the MP algorithms converge much faster; the empirical
complexity for computing the equilibrium flows up to a certain error
$|\boldsymbol{x}^{MP}-\boldsymbol{x}^{*}|$ is roughly $O(|E|)$.
However, there exists some networks where MP with Method II does not
converge, while MP with Method I converges successfully. For the experiments
in the main text, we use Method I to treat the destination node for its
better convergence properties. 

\subsection{Extension to The Case of Multiple Destination\label{subsec:lower_MP_multiclass}}

The case of multiple destinations can be studied similarly. The traffic
flows can be classified into different classes according to their
destinations. Let $N_{d}$ denotes the number of destinations, and
$x_{e}^{a}$ denote the flow on edge $e$ targeted at the $a$-th
destination (or the $a$-th class), the lower-level optimization problem
(for solving equilibrium flows) is defined as
\begin{align}
\min_{\boldsymbol{x}} \;\; &  \Phi(\boldsymbol{x}) =\sum_{e\in E}\int_{0}^{\sum_{a}x_{e}^{a}}\ell_{e}(y)\mathrm{d}y=\sum_{e\in E}\phi_{e}\bigg(\sum_{a=1}^{N_{d}}x_{e}^{a}\bigg), \\
& \text{s.t. }R_{i}^{a} :=\Lambda_{i}^{a}+\sum_{e\in\partial i}B_{i,e}x_{e}^{a}=0,\quad\forall i,a, \\
& \qquad \qquad x_{e}\geq0,\quad\forall e,a.
\end{align}

\begin{figure}
\includegraphics[scale=1.5]{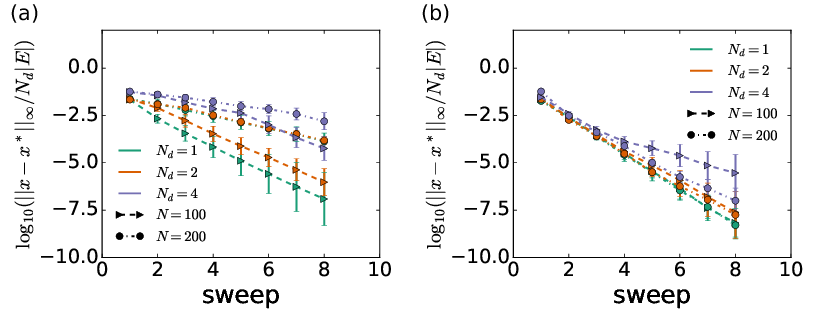}

\caption{Message-passing algorithm for routing games with $N_{d}$ destinations
converges to the equilibrium flows $\boldsymbol{x}^{*}$. (a) Method
I in Sec.~\ref{subsec:treatmen_of_destination}) is used to treat
the destination nodes. (b) Method II in Sec.~\ref{subsec:treatmen_of_destination})
is used to treat the destination nodes. Each sweep comprises of $40|E|$
local MP updates. \label{fig:routing_game_Lower_MP_converge_multi_classes}}
\end{figure}

To accommodate the nonlinear interactions of flows of different classes
$\{x_{e}^{a}\}$ in $\phi_{e}(\sum_{a}x_{e}^{a})$, we adopt a coordinate-descent
like approach in the MP algorithm as follows. In the treatment of
the flows of class $a$, we fix the flows of other classes $\{x_{e}^{b}\}_{b\neq a}$
to their working points and compute the messages of class $a$ as
\begin{align}
\Phi_{i\to e}^{a}(x_{e}^{a}) & =\min_{\{x_{e'}^{a}\ge0\}|R_{i}^{a}=0}\sum_{e'\in\partial i\backslash e}\bigg[\Phi_{k\to e'}^{a}(x_{e'}^{a})+\phi_{e'}\big(\sum_{b\neq a}\tilde{x}_{k\to e'}^{b}+x_{e'}^{a}\big)\bigg]\nonumber \\
 & =\min_{\{x_{e'}^{a}\ge0\}|R_{i}^{a}=0}\sum_{e'\in\partial i\backslash e}\bigg[\Phi_{k\to e'}^{a}(\tilde{x}_{k\to e'}^{a}+\varepsilon_{e'}^{a})+\phi_{e'}\big(\sum_{b}\tilde{x}_{k\to e'}^{b}+\varepsilon_{e'}^{a}\big)\bigg].
\end{align}
We further adopt the approximations $\Phi_{k\to e'}^{a}(\tilde{x}_{k\to e'}^{a}+\varepsilon_{e'}^{a})\approx\text{\ensuremath{\Phi_{k\to e'}^{a}}(\ensuremath{\tilde{x}_{k\to e'}^{a}}})+\beta_{k\to e'}^{a}\varepsilon_{e'}^{a}+\frac{1}{2}\alpha_{k\to e'}^{a}(\varepsilon_{e'}^{a})^{2}$
(augmented by a piecewise quadratic function if it is non-smooth)
and $\phi_{e'}(\sum_{b}\tilde{x}_{k\to e'}^{b}+\varepsilon_{e'}^{a})\approx\phi_{e'}(\sum_{b}\tilde{x}_{k\to e'}^{b})+\phi_{e'}^{\prime}(\sum_{b}\tilde{x}_{k\to e'}^{b})\varepsilon_{e'}^{a}+\frac{1}{2}\phi_{e'}^{\prime\prime}(\sum_{b}\tilde{x}_{k\to e'}^{b})(\varepsilon_{e'}^{a})^{2}$,
and solve for the coefficients $\{\alpha_{i\to e}^{a},\beta_{i\to e}^{a}\}$
as in the single-destination scenario. The resulting MP algorithm
has the same structure as the one of single class described above,
where its efficacy is shown in Fig.~\ref{fig:routing_game_Lower_MP_converge_multi_classes}. 

\subsection{Bilevel Optimization in Routing Games}

The toll optimization problem (for single destination) of the upper-level planner is defined as
\begin{equation}
\min_{\boldsymbol{\tau}}H(\boldsymbol{x}^{*}(\boldsymbol{\tau}))=\sum_{e\in E}x_{e}^{*}(\boldsymbol{\tau})\ell_{e}\big(x_{e}^{*}(\boldsymbol{\tau})\big),
\end{equation}
\begin{equation}
\text{s. t. constraints of }\boldsymbol{\tau} \text{ and:}
\end{equation}
\begin{align}
\boldsymbol{x}^{*}(\boldsymbol{\tau})= & \text{arg}\min_{\boldsymbol{x}} \Phi(\boldsymbol{x}; \boldsymbol{\tau}) =  \text{arg}\min_{\boldsymbol{x}}\sum_{e\in E}\int_{0}^{x_{e}}\big[\ell_{e}(y)+\tau_{e}\big]\mathrm{d}y,\label{eq:xN_of_tau_lower_level_def}\\
 & \text{s.t. }x_{e}\geq0,R_{i}=0,\forall e,i.
\end{align}
Here, the social cost $H(\boldsymbol{x})$ is the objective function of the upper-level problem, i.e., the overall target of the central planner is to reduce this social cost. The potential function $\Phi(\boldsymbol{x}; \boldsymbol{\tau})$ is the objective function of the lower-level problem, which governs the equilibrium flow $\boldsymbol{x}^{*}(\boldsymbol{\tau})$ for a given toll configuration $\boldsymbol{\tau}$.

The computation of the social optimum $H(x)$ has a similar form for computing the equilibrium flows, we therefore use a parallel MP procedure for the social cost
\begin{align}
H_{i\to e}(x_{e}) & =\min_{\{x_{e'}\geq0\}|R_{i}}\sum_{e'\in\partial i\backslash e}\bigg\{ H_{k\to e'}(x_{e'})+\sigma(x_{e'})\bigg\},\nonumber \\
 & =\min_{\{x_{e'}\geq0\}|R_{i}}\sum_{e'\in\partial i\backslash e}\bigg\{\frac{1}{2}\gamma_{k\to e'}\big(\varepsilon_{e'}\big)^{2}+\delta_{k\to e'}\varepsilon_{e'}+\frac{1}{2}\sigma_{e'}''(\tilde{x}_{k\to e'})\big(\varepsilon_{e'}\big)^{2}+\sigma_{e'}'(\tilde{x}_{k\to e'})\varepsilon_{e'}\bigg\},
\end{align}
where $H_{k\to e'}(x_{e'})$ assumes a quadratic approximation and
needs to be augmented by a piecewise quadratic function if it is non-smooth.
It results in an upper-level MP algorithm having the same structure
as the one for computing the equilibrium flows. The difference is that,
since the flow is not directly driven by the central planner, the
working point $\tilde{x}_{i\to e}$ is not updated at the upper level
MP. Instead, the central planner updates the toll $\tau_{e}$
such that selfish users are attracted to the solution with a lower
social cost.

When the toll $\tau_{e}$ on edge $e$ is adapted, the marginal Nash-equilibrium
flow $x_{e}^{N}$ changes accordingly
\begin{equation}
x_{e}^{N}(\tau_{e})=\text{arg}\min_{x_{e}\ge0}\big[\Phi_{i\to e}(x_{e})+\Phi_{j\to e}(x_{e})+\phi_{e}(x_{e})+\tau_{e}x_{e}\big].
\end{equation}
For smooth message functions $\Phi_{i\to e}(x_{e})$ and $\Phi_{j\to e}(x_{e})$
\begin{equation}
x_{e}^{N}(\tau_{e})=\max\bigg(\frac{\big(\alpha_{i\to e}+\frac{1}{2}\phi''_{i\to e}\big)\tilde{x}_{i\to e}+\big(\alpha_{j\to e}+\frac{1}{2}\phi''_{j\to e}\big)\tilde{x}_{j\to e}-(\beta_{i\to e}+\beta_{j\to e}+\tau_{e}+\frac{1}{2}\phi'_{i\to e}+\frac{1}{2}\phi'_{j\to e})}{\alpha_{i\to e}+\alpha_{j\to e}+\frac{1}{2}\phi''_{i\to e}+\frac{1}{2}\phi''_{j\to e}},0\bigg),
\end{equation}
which is a piecewise linear function of $\tau_{e}$ with two branches.
For non-smooth cavity functions, $x_{e}^{N}(\tau_{e})$ can also be
obtained straightforwardly, which is a piecewise linear function of
$\tau_{e}$ with multiple branches.

The goal of toll-adaptation of $\tau_{e}$ is to decrease the social
cost $H(\boldsymbol{x})$, which amounts to decrease the full social
cost on edge $e$
\begin{align}
\tau_{e}^{*} & =\text{arg}\min_{\tau_{e}}H_{e}^{\text{full}}(x_{e}^{N}(\tau_{e})),\label{eq:opt_He_full}\\
H_{e}^{\text{full}}(x_{e}) & :=H_{i\to e}(x_{e})+H_{j\to e}(x_{e})+\sigma_{e}(x_{e}),
\end{align}
where the optimization in Eq.~(\ref{eq:opt_He_full}) needs to obey
necessary constraints on tolls (e.g., the restriction $0\leq\tau_{e}\leq\tau_{e}^{\max}$
is considered in the main text).

As $H_{e}^{\text{full}}(x_{e}^{N}(\tau_{e}))$ is a convex function
of $x_{e}^{N}$, it is sufficient to adapt $\tau_{e}$ such that $x_{e}^{N}(\tau_{e})$
gets as close to the marginal socially optimal flow $x_{e}^{G}$ as
possible, where $x_{e}^{G}$ is given by
\begin{align}
x_{e}^{G} & =\text{arg}\min_{x_{e}\ge0}\big[H_{i\to e}(x_{e})+H_{j\to e}(x_{e})+\sigma_{e}(x_{e})\big].
\end{align}
The search for the optimal toll $\tau_{e}^{*}$ can be done efficiently
by utilizing the property that $x_{e}^{N}(\tau_{e})$ is a piecewise
linear function of $\tau_{e}$.

The resulting bilevel MP algorithm is described in the
main text, where the lower-level messages $\{\alpha_{i\to e}^{(m)},\beta_{i\to e}^{(m)},\tilde{x}_{i\to e}\}$ ($m\in\{L,R\}$) and upper-level messages $\{\gamma_{i\to e}^{(m)},\delta_{i\to e}^{(m)}\}$
are passed along edges to compute the equilibrium flows $x_{e}^{N}$
and related quantities. These messages facilitate the computation
of $H_{e}^{\text{full}}(x_{e})$ in the upper level, which is used
to update the toll variables $\tau_{e}$. In practice, the update of tolls is
less frequent then the update of other messages. In the experiments
shown in the main text, for every $\frac{2}{5}N_{d}|E|$ MP iterations,
we randomly select an edge $e$ and update its toll.

The resulting bilevel MP algorithm is summarized in Algorithm~\ref{alg:MP_routing_bilevel}.

\begin{algorithm}
\SetNoFillComment
\caption{Message-passing algorithm for toll optimization in routing games (bilevel, single destination, Method I to treat the destination node $\mathcal{D}$ )}\label{alg:MP_routing_bilevel}

\KwIn{Road network $G(V,E)$ (pre-processed to remove dangling nodes), node parameters $\{ \Lambda_i \}$ defining the resources, edge parameters defining $\ell_e(x_e)$, $\phi_e(x_e)$ and $\sigma_e(x_e)$, maximal number of iterations $T$, time interval $t_{\text{update\_intv}}$ for updating tolls, time interval $t_{\text{dump\_intv}}$ for saving intermediate tolls.}

\BlankLine
Initialize the messages $\{ \alpha^{(m)}_{i \to e}, \beta^{(m)}_{i \to e} \}$ and $\{ f^{\text{lw}}_{i \to e} \}$ (effective leaf states of the lower layer) randomly.

\BlankLine
Initialize the messages $\{ \gamma^{(m)}_{i \to e}, \delta^{(m)}_{i \to e} \}$ and $\{ f^{\text{up}}_{i \to e} \}$ (effective leaf states of the upper layer) randomly.

\BlankLine
Initialize the working points $\{ \tilde{x}_{i \to e} \}$ and tolls $\{ \tau_{e} \}$.

\BlankLine
\For{$t \textnormal{ in } 1:T$}{

Randomly select a node $i$ and one of its adjacent edge $e \in \partial i$.

\BlankLine
\BlankLine
\tcc{update the lower-level messages:}
Run Subroutine (a) in Algorithm~\ref{alg:MP_routing_lower} to update the messages $\alpha^{(m)}_{i \to e}, \beta^{(m)}_{i \to e}, f^{\text{lw}}_{i \to e}$ using $\phi_e(\cdot)$.

\BlankLine
\If{$f_{i \to e} = 1$}{
Set $\tilde{x}_{i \to e} = \Lambda^{\text{eff}}_{i \to e}$.
}

\BlankLine
\BlankLine
\tcc{update the upper-level messages:}

\uIf{$i = \mathcal{D}$}{
Set $\gamma_{i \to e} = 0, \delta_{i \to e} = 0, f^{\text{up}}_{i \to e} = 0$.
}

\Else
{Replace $\{ \alpha^{(m)}_{i \to e}, \beta^{(m)}_{i \to e}, f^{\text{lw}}_{i \to e}$ using $\phi_e(\cdot) \}$ by $\{ \gamma^{(m)}_{i \to e}, \delta^{(m)}_{i \to e}, f^{\text{up}}_{i \to e} \}$ in Sec.~\ref{sec:MP_for_nonsmooth}. 

Follow the same procedures therein to update the messages $\gamma^{(m)}_{i \to e}, \delta^{(m)}_{i \to e}, f^{\text{up}}_{i \to e}$ with the working point evaluated at the up-to-date value of $\tilde{x}_{i \to e}$.
}

\BlankLine
\BlankLine

\tcc{update the working points:}
Using the up-to-date messages $\{ \alpha^{(m)}_{i \to e}, \beta^{(m)}_{i \to e}, f^{\text{lw}}_{i \to e} \}$, $\{ \tilde{x}_{i \to e} \}$ and $\{ \tau_e \}$ to determine the approximated forms of $\Phi_{i \to e}(x_e)$ and $\Phi_{j \to e}(x_e)$.

Compute $x_e^*=\arg \min _{x_e \geq 0}\left[\Phi_{i \to e}\left(x_e\right)+\Phi_{j \to e}\left(x_e\right)+\phi_e\left(x_e\right)\right]$.

Update the working point as $\tilde{x}_{i \to e}^{\text {new }} \leftarrow s x_e^*+(1-s) \tilde{x}_{i \to e}^{\text {old }}$ with a learning rate $s$.

\BlankLine
\BlankLine
\tcc{update the tolls:}
\If{$t \; \operatorname{mod} \; t_{\textnormal{update\_intv}} = 0$}
{Using the up-to-date messages $\{ \gamma^{(m)}_{i \to e}, \delta^{(m)}_{i \to e}, f^{\text{up}}_{i \to e} \}$ and $\{ \tilde{x}_{i \to e} \}$ to determine the approximated forms of $H_{i \to e}(x_e)$ and $H_{j \to e}(x_e)$.

Compute $x_e^G = \arg \min _{x_e \geq 0}\left[H_{i \to e}\left(x_e\right)+ H_{j \to e}\left(x_e\right)+ \sigma_e\left(x_e\right)\right]$.

Determine the (estimated) equilibrium flow as a function of toll: $x_e^N(\tau_e)=\arg \min _{x_e \geq 0}\left[\Phi_{i \to e}\left(x_e\right)+\Phi_{j \to e}\left(x_e\right)+\phi_e\left(x_e\right | \tau_e)\right]$, which is a piecewise linear function.

Find $\tau_e \in [0, \tau_e^{\max}]$ such that $x_e^*(\tau_e)$ is as close to $x_e^G$ as possible, and update $\tau_e$ accordingly.
}

\BlankLine
\BlankLine
\If{$t \; \operatorname{mod} \; t_{\textnormal{save\_intv}} = 0$}
{Save the up-to-date tolls.}

\BlankLine
\If{\textnormal{messages converge}}{
Exit the for loop.
}

} 

\BlankLine

For each edge $e$, compute the equilibrium flow as $x_e^*=\arg \min _{x_e \geq 0}\left[\Phi_{i \to e}\left(x_e\right)+\Phi_{j \to e}\left(x_e\right)+\phi_e\left(x_e\right)\right]$.

\BlankLine
\KwOut{Convergence status, equilibrium flows $\{ x^{*}_{e} \}$, the optimized tolls $\{ \tau_e \}$}
\end{algorithm}

\subsubsection{Extension to Multiple Destinations}

\begin{figure}
\includegraphics[scale=1.5]{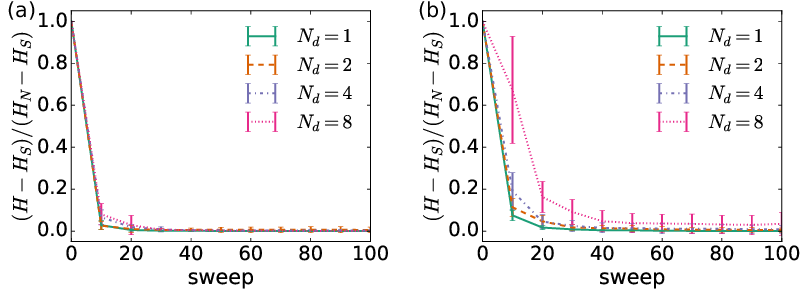}

\caption{The effect of tolls on the reduction in fractional social cost in routing
games on random regular graphs with $N_{d}$ destinations, where tolls
are restricted as $0\protect\leq\tau_{e}\protect\leq1$. Each sweep
comprises $40N_{d}|E|$ local MP updates and $100$ edgewise
toll updates in a random sequential schedule. (a) $N=100$. (b) $N=200$.
\label{fig:results_toll_opt_multiclass}}
\end{figure}

The toll optimization problems with multiple destinations can be tackled
by the proposed bilevel MP algorithm using the approximations
in Sec.~\ref{subsec:lower_MP_multiclass}. The results are shown in
Fig.~\ref{fig:results_toll_opt_multiclass}, which demonstrate the
effectiveness of the algorithm in reducing the social cost by adapting
tolls.

\subsection{The Bilevel Programming Approach}

Here we demonstrate the results of the bilevel programming approach
to the toll optimization problem. It is achieved by expressing the solution
of Eq.~(\ref{eq:xN_of_tau_lower_level_def}) as constraints imposed
by the Karush--Kuhn--Tucker (KKT) condition, and solve the bilevel
optimization as a global nonlinear programming problem~\cite{Colson2007}.
Such an approach is intrinsically difficult as (i) the constraints
by the KKT conditions can be nonlinear and non-convex; (ii) the complementary
slackness conditions are combinatorial, which requires a treatment
with mixed integer programming (e.g., through branch and bound).
Therefore, the bilevel programming approach offers a centralized algorithm that is generally non-scalable.

\begin{figure}
\includegraphics[scale=0.4]{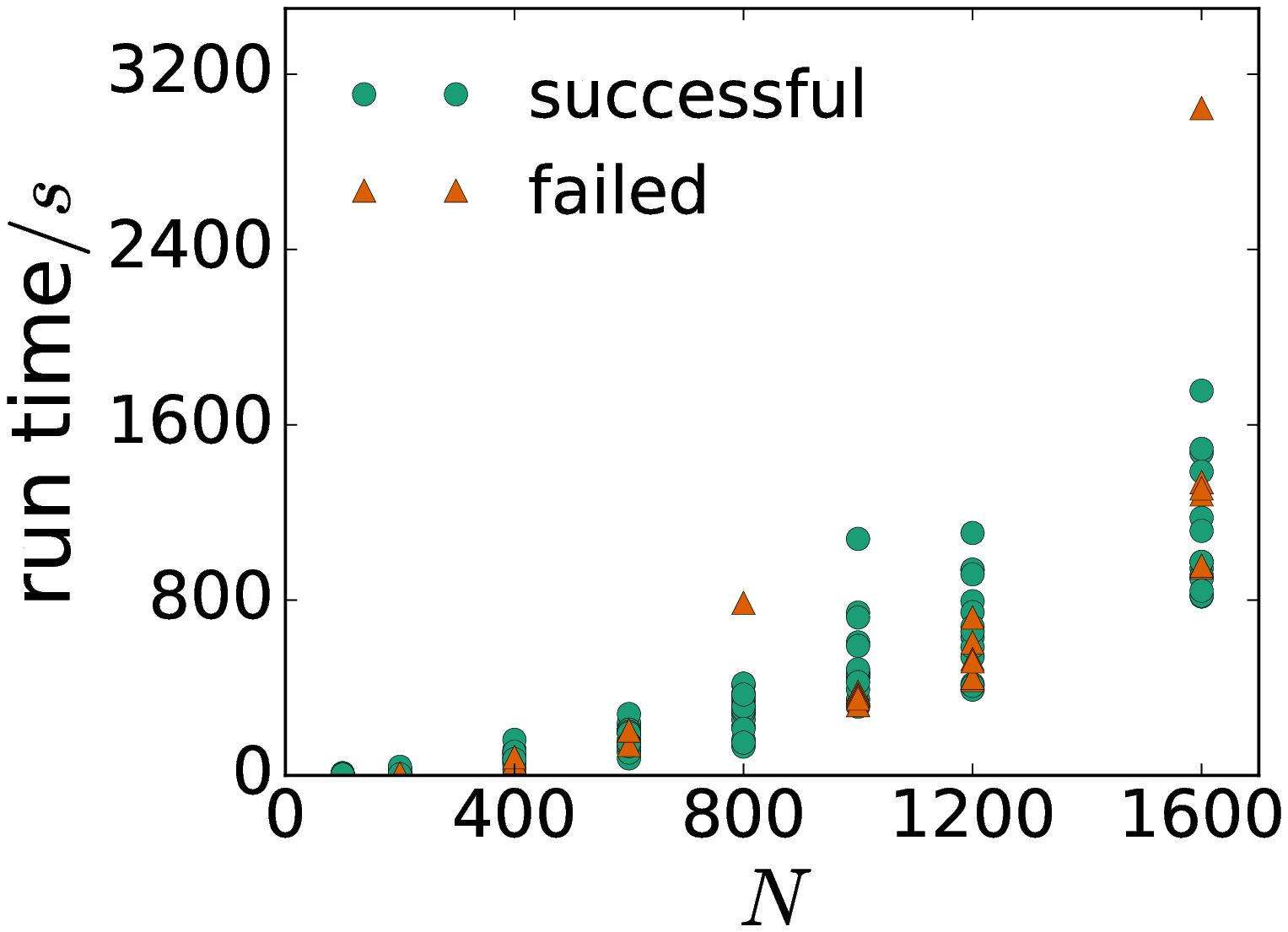}

\caption{Run time of bilevel programming on the toll optimization problem in
random regular graphs. For each network size, 20 different problem
realizations are considered; red triangles represent the cases
where bilevel programming fails to find the solution in single trials, green dots represent successful trials.\label{fig:runtime_BLP}}

\end{figure}

It is difficult to directly compare the MP algorithms to the bilevel
programming approach, as the convergence rate for either approach
is difficult to establish. Besides, the bilevel programming approach
is a centralized optimization method, which has a different space complexity
per iteration. Nevertheless, we present the results of CPU run time
(CPU in used: i5-3317U) of bilevel programming (package in used: bileveljump.jl~\cite{Joaquim2021} with the IPOPT solver~\cite{Wachter2006}) on the toll optimization problem in Fig.~\ref{fig:runtime_BLP}.
There exist cases where bilevel programming fails to find the solution
in a single trial, and the run times vary significantly among different
problem realizations.

The bilevel programming approach is more generic and flexible than
the MP approach, but it does not offer a decentralized
algorithm as the MP approach. Besides, the MP algorithms can be extended
to the scenarios with discrete variables, which is very difficult
for the global optimization approach. It is also difficult to treat
the toll selection problem in the main text with the bilevel programming
method, especially when the socially optimum is not known a priori
in some variants of toll-setting or network-design problems~\cite{Migdalas1995}.

\section{Message-passing Algorithms For Atomic Routing Games\label{sec:MP_atomic}}

Message-passing algorithms can be extended to include atomic routing games, where each player controls one unit of traffic. The atomic games differ from the above-studied non-atomic games in that the flow variable of edge $e$ is an integer $x_{e} \in \mathbb{Z}$. 
The integer constraints make even the single-level optimization a difficult integer programming problem. An existing mixed integer programming approach can be used to solve the equilibrium flow problem of atomic routing games for moderate network sizes. 
However, the bilevel optimization (such as the toll-setting problem) appears to be much more difficult to solve.

Message-passing algorithms have also been successfully applied to solve network flow and routing problems (as a single-level optimization problem) with integer flow variables~\cite{Yeung2012,Yeung2013IEEE,Bacco2014,Yeung2014}. Here, we generalize the MP algorithms for integer flows to solve the Wardrop equilibrium problem in directed networks.
The problem of atomic games is defined similarly as the non-atomic games introduced in Sec.~\ref{subsec:def_nonatomic}, except that there are additional integer constraints $x_{e} \in \mathbb{Z}$, and that the potential function $\Phi(\boldsymbol{x})$ (to be minimized) in Eq.~(\ref{eq:Phi_def_nonatomic}) is replaced by
\begin{equation}
\Phi(\boldsymbol{x}) = \sum_{e\in E} \sum_{y=1}^{x_{e}} \big[ \ell_{e}(y) +\tau_{e} \big] =:\sum_{e\in E}\phi_{e}(x_{e}).\label{eq:Phi_def_atomic}
\end{equation}

The MP equation for minimizing the potential $\Phi(\boldsymbol{x})$ can be written similarly to the non-atomic games as
\begin{align}
\Phi_{i\to e}(x_{e}) & =\min_{\{ x_{e'} \in \mathbb{N} \}|R_{i}=0}\sum_{e'\in\partial i\backslash e}\bigg[\Phi_{e'\to i}(x_{e'})+\phi_{e'}(x_{e'})\bigg], \label{eq:Phi_cav_atomic}
\end{align}
where the message function $\Phi_{i \to e}(x_{e})$ is defined on a 1-dimensional grid $x_{e} \in \{0,1,2,...\}$ to be computed recursively.

To simplify the above message-passing calculations, we can adopt the perturbation approach similar to the case of non-atomic game, by considering only a few flow values $x_{e}$ near a certain working point $\tilde{x}_{i \to e} \in \mathbb{N}$~\cite{Yeung2013IEEE}
\begin{equation}
\Phi_{i \to e}(x_e) = \Phi_{i \to e}(\tilde{x}_{i \to e} + m) =: h^{m}_{i \to e}, \qquad m\in\mathbb{Z}, -M \leq m \leq M,
\end{equation}
where $M$ is a small integer parameter determining the scope to look around $\tilde{x}_{i \to e}$. The MP algorithm proceeds by iteratively solving Eq.~(\ref{eq:Phi_cav_atomic}) for different directed edges $\{ i \to e | i \in V, e \in E \}$. The optimization problem for a particular edge $i \to e$ is approximately achieved by searching the grid $\{ x_{e'} \in \mathbb{N} \}_{\forall e' \in \partial i \backslash e}$ in the vicinity of $\{ \tilde{x}_{e' \to i} \}_{\forall e' \in \partial i \backslash e}$. The complexity of each iteration is proportional to $(2M+1)^{|\partial i|-1}$. In our experiments, we use $M=1$ following the suggestion of \cite{Yeung2013IEEE}. To some extent, the implementation is more straightforward than that in non-atomic games, as there is no need to treat the non-smooth cavity energy function.

Similarly to the non-atomic game case, the working point $\tilde{x}_{i \to e}$ is also gradually pushed towards the minimizer of the edgewise full energy 
\begin{align}
x_{e}^{*} & = \text{arg}\min_{x_{e} \in \mathbb{N}} \Phi^{\text{full}}_{e}(x_{e})  = \text{arg}\min_{x_{e} \in \mathbb{N}}\big[\Phi_{i\to e}(x_{e})+\Phi_{j\to e}(x_{e})+\phi_{e}(x_{e})\big], \label{eq:xe_marginal_def_atomic} \\
\tilde{x}_{i \to e}^{\text{new}} & \leftarrow \tilde{x}_{i \to e}^{\text{old}} + \text{sign}(x^{*}_{e} - \tilde{x}_{i \to e}^{\text{old}}), \label{eq:update_working_point_atomic}
\end{align}
in which case the working point $\tilde{x}_{i \to e}$ is expected to get closer and closer to the equilibrium point $x^{*}_{e}$ and the perturbation scheme will become more accurate.

In Eq.~(\ref{eq:xe_marginal_def_atomic}), we can only access the flow values in the set $A_{e} := \{ \tilde{x}_{i\to e} +m | -M \leq m \leq M \} \cap \{ \tilde{x}_{j\to e} +m | -M \leq m \leq M \} $. If $A_{e} = \varnothing$, instead of applying Eq.~(\ref{eq:update_working_point_atomic}), we push $\tilde{x}_{i \to e}$ towards $\tilde{x}_{j \to e}$ incrementally to increase the chance of overlap in future steps
\begin{equation}
\tilde{x}_{i \to e}^{\text{new}} \leftarrow \tilde{x}_{i \to e}^{\text{old}} + \text{sign}(\tilde{x}_{j \to e}^{\text{old}} - \tilde{x}_{i \to e}^{\text{old}}).
\end{equation}

The message-passing algorithm in our study is known as the min-sum algorithm, where caution is needed when there are degenerate energy minima such that taking the minimum is ambiguous (see Sec.~8.4.5 of \cite{Bishop2006} for example). This issue does not impact on MP algorithms for non-atomic games due to the convex nature of the single-level problem, but complicates the use of MP algorithms for atomic games. In~\cite{Yeung2012, Yeung2013IEEE}, a small random bias field $\xi_{e}$ is added to non-linear cost $\phi_{e}(x_{e})$ to break the degeneracy
\begin{equation}
\phi_{e}'(x_{e}) \leftarrow \phi_{e}(x_{e}) + \xi_{e} |x_{e}|.
\end{equation}

Alternatively, we can also select one of the degenerate solution based on its consistency with the constraints. For example, in Eq.~(\ref{eq:xe_marginal_def_atomic}), suppose there are two flow values $\{ x^{*(1)}_{e}, x^{*(2)}_{e} \}$ corresponding to $\min \Phi^{\text{full}}_{e}(x_{e})$. We would choose the solution that is more consistent with the flow conservation constraint
\begin{align}
R_{i \to e}^{(n)} & = \Lambda_{i} + \sum_{e'=(k,i) \in \partial i \backslash e} B_{i,e'} \tilde{x}_{k \to e'} + B_{i,e} x_{e}^{*(n)}, \quad n = 1,2, 
\end{align}
that is, we would assign the solution $x^{*(n)}_{e}$ with minimal value of $|R_{i \to e}^{(n)}|$ to $x_{e}^{*}$ in Eq.~(\ref{eq:xe_marginal_def_atomic}).

Both methods facilitate the convergence of the MP algorithms to a local minima in atomic routing games. We refer readers to Ref.~\cite{Yeung2013IEEE} for more details on MP algorithms in network flow problems with integer constraints.

\subsection{Bilevel MP Algorithms for Atomic Games}

The method to treat bilevel optimization (in particular, the toll-setting problem) in non-atomic games can also be applied to atomic games. This is achieved by considering a parallel message-passing process for minimizing the social cost (in the upper-level)
\begin{align}
H_{i\to e}(x_{e}) & =\min_{\{ x_{e'} \in \mathbb{N} \}|R_{i}=0}\sum_{e'\in\partial i\backslash e}\bigg[H_{e'\to i}(x_{e'})+\sigma_{e'}(x_{e'})\bigg], \label{eq:H_cav_atomic}
\end{align}
where $\sigma_{e'}(x_{e'}) = x_{e'} \ell_{e}(x_{e'})$.
The minimizer of the edgewise full cost $x^{S}_{e} = \text{arg}\min_{x_{e}} H^{\text{full}}_{e}(x_{e}) = \text{arg}\min_{x_{e}} \big[ H_{i\to e}(x_{e}) + H_{j\to e}(x_{e}) + \sigma_{e}(x_{e}) \big]$ is informative of the min social-cost flow during the upper-level MP updates. We then update the toll $\tau_{e}$ incrementally such that the toll-dependent equilibrium flow $x_{e}^{*}(\tau)$ (given by Eq.~(\ref{eq:xe_marginal_def_atomic})) gets closer to $x^{S}_{e}$.

\begin{figure}
    \centering
    \includegraphics[scale=0.4]{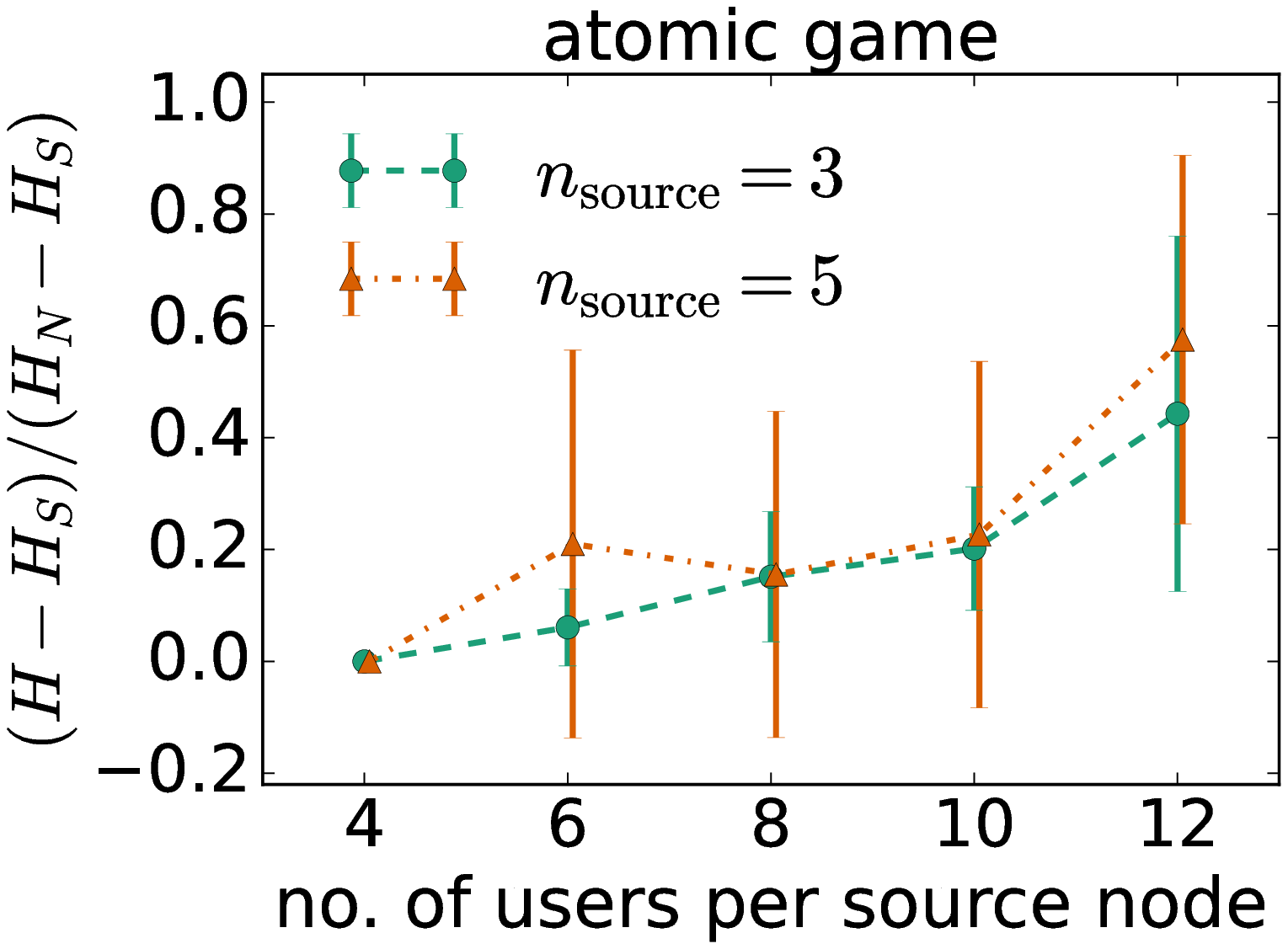}
    \caption{The effect of tolls on the fractional reduction in social cost. Random regular graphs of degree 3 and size $N=50$ are considered. Each data point is the average of 10 different problem realizations. For each problem realization, we consider 5 different trials of bileve MP processes (with different random starting points), and recorded a few of the tolls during each process; we then kept the toll configuration corresponding to the smallest social cost. In each experiment, three or five nodes (dashed lines or dash-dotted lines) are randomly selected as the sources (origins), on which the users aim to route to a universal destination.}
    \label{fig:dH_MP_atomic}
\end{figure}

We demonstrate the effectiveness of the bilevel MP algorithm for toll-setting problems in atomic games in Fig.~\ref{fig:dH_MP_atomic}. We observe that the algorithm is effective for small networks with light and moderate loads, while it becomes less effective for heavy loads. The performance also deteriorate in large networks (not shown). We conjecture that a more non-local toll update method is needed to improve the performance in such cases, which remains to be explored in future studies. Nevertheless, the bilevel optimization problems in atomic games are intrinsically difficult combinatorial optimization problems, where the complex energy landscape is prohibitive for most optimization algorithms.

\subsection{Optimal Toll Configurations for Some Instances}

In this section, we demonstrate the optimal tolls found by the bilevel MP algorithms for some cases in the Sioux Falls road network~\cite{Suwansirikul1987}, which is a popular benchmark network for transportation research, consisting of 24 nodes and 76 directed edges.

In this example, upon obtaining the optimal tolls $\{ \tau_{e}^{\text{opt}} \}$ in the bilevel MP algorithm, we keep only the tolls with relatively large values $\tau_{e}^{\text{new}} \leftarrow \tau_{e}^{\text{opt}} \Theta(\tau_{e}^{\text{opt}} \geq \epsilon)$ (where $\epsilon$ is a small threshold). If this does not incur an increment in social cost, we keep $\{ \tau_{e}^{\text{new}} \}$, otherwise, we keep $\{ \tau_{e}^{\text{opt}} \}$.

\begin{figure}
    \centering
    \includegraphics[scale=1.6]{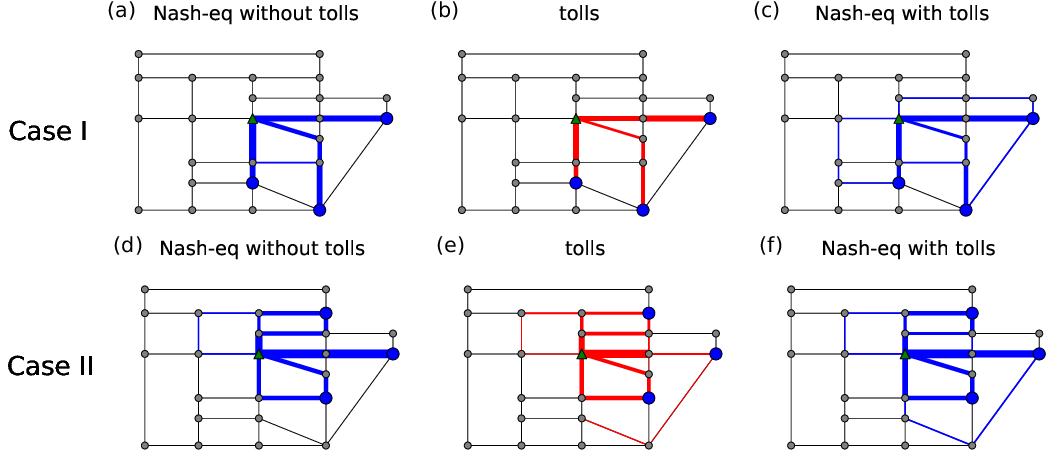}
    \caption{Optimal tolls and network flows in the Sioux Falls network found by the bilevel MP algorithm. Note that the underlying network is a directed graph (each undirected edge in the figures comprises two directed edges in both directions). We set the central node as the destination node (marked as a green triangle), and randomly select 3 nodes as the source nodes (marked as blue circles) in either case. In Case I (top), each source node has 4 users, while in Case II (bottom), each source has 6 users. In (a)(c) and (d)(f), the network flow patterns in Nash (Wardrop) equilibrium are shown, where a non-zero flow is marked as a blue edge (the edge width is proportional to the flow magnitude). In (b) and (e), the optimal tolls are shown, where a non-zero toll is marked as a red edge (the edge width is proportional to the toll magnitude).}
    \label{fig:SiouxFalls_case}
\end{figure}

The results are shown in Fig.~\ref{fig:SiouxFalls_case}. In Case I, we observed that the shortest paths from the source nodes to destination are very congested, and the optimal tolls are placed on these paths such that some users are re-routed to alleviate congestion. In Case II, interestingly, tolls are also placed on some non-shortest paths in the optimal configuration found by the bilevel MP algorithm.

\section{Message-Passing Algorithms for Flow Control in Undirected Networks}

In this section, we provide the details of the MP algorithm
for flow control in undirected networks. In a simple undirected graph
$G(V,E)$, nodes $i$ and $j$ can be connected by at most one edge
$(i,j)$, where the order of node $i$ and $j$ does not matter (in contrast, edges
$(i,j)$ and $(j,i)$ are two different edges in a directed
graph). In this case, edge $(i,j)$ can either transmit resources
from node $j$ to node $i$ or from node $i$ to node $j$. Denoting
$x_{ij}$(=$-x_{ji}$) as the flow from node $j$ to node $i$; if
$x_{ij}<0$ (or $x_{ji}=-x_{ij}>0$), the resources are being
transmitted from node $i$ to node $j$. We also assume that the underlying
graph does not have any leaf nodes, by recursively trimming leaf
nodes and absorbing their resources into neighboring nodes.

\subsection{MP For Lower-level Optimization}
\label{sec:MPLowLevelOpt}
The equilibrium flow (in the lower-level optimization problem) is
the minimizer of the problem
\begin{align}
\min_{\boldsymbol{x}}C(\boldsymbol{x})= & \sum_{(i,j)}\frac{1}{2}r_{ij}x_{ij}^{2},\label{eq:undirected_lower_level_prob}\\
\text{s.t. }R_{i}= & \Lambda_{i}+\sum_{j\in\mathcal{N}_{i}}x_{ij}=0,\quad\forall i\neq\mathcal{D},\label{eq:undirected_flow_conservation}
\end{align}
where the reference node $\mathcal{D}$ can be arbitrarily chosen.

The above optimization problem can be mapped onto its dual problem as
\begin{align}
\min_{\boldsymbol{\mu}}C^{\text{dual}}(\boldsymbol{\mu}) & =\sum_{(i,j)}\frac{1}{2r_{ij}}(\mu_{j}-\mu_{i})^{2}-\sum_{i}\Lambda_{i}\mu_{i},\\
 & =:\frac{1}{2}\boldsymbol{\mu}^{\top}L\boldsymbol{\mu}-\boldsymbol{\Lambda}^{\top}\boldsymbol{\mu},
\end{align}
where $\mu_{i}$ is the Lagrange multiplier (or dual variable)
associated with the flow conservation constraint $R_{i}=0$, and $L$
is the Laplacian matrix with matrix element
\begin{equation}
L_{ij}:=\left(\sum_{k\in\mathcal{N}_{i}}\frac{1}{r_{ik}}\right)\delta_{ij}-\frac{1}{r_{ij}}.\label{eq:lapl_matrix_def}
\end{equation}

The solution of the dual problem can be obtained by solving the system
of linear equations $L\boldsymbol{\mu}^{*}=\boldsymbol{\Lambda}$,
and the equilibrium flow $x_{ij}^{*}$ is related to the optimal Lagrange
multiplier $\boldsymbol{\mu}^{*}$ through $x_{ij}^{*}=\frac{\mu_{j}^{*}-\mu_{i}^{*}}{r_{ij}}$.
The drawback of such an approach is that (i) solving the systems of
linear equations usually needs a centralized solver; (ii) to compute
the response of the equilibrium flow to changes of the control parameters
$\{r_{ij}\}$ in bilevel optimization, one needs to evaluate the pseudo-inverse
of the Laplacian matrix per iteration, which can be very computationally
demanding for large networks. Instead, we proposed to use MP
for computing the equilibrium flows and tackle the related
bilevel optimization problem, which is a scalable and efficient
decentralized algorithm.

For the lower-level equilibrium flow problem in Eq.~(\ref{eq:undirected_lower_level_prob}),
the MP algorithm amounts to computing the message functions
\begin{equation}
C_{i\to j}(x_{ij})=\min_{\{x_{ki}\}|R_{i}=0}\bigg[\frac{1}{2}r_{ij}x_{ij}^{2}+\sum_{k\in\mathcal{N}_{i}\backslash j}C_{k\to i}(x_{ki})\bigg],\label{eq:cavity_C_ij}
\end{equation}
where the definition of the message function $C_{i\to j}(x_{ij})$
differs from the one in Eq.~(\ref{eq:cavity_Phi_ie}) in that it includes
the interaction term on edge $(i,j)$, which yields a more concise
update rule in this problem. Similar to routing games, we approximate
the message function by a quadratic form $C_{i\to j}(x_{ij})=\frac{1}{2}\alpha_{i\to j}(x_{ij}-\hat{x}_{i\to j})^{2}+\text{const}$,
such that the local optimization in Eq.~(\ref{eq:cavity_C_ij}) reduces
to the computation of the real-number messages $m_{i\to j}\in\{\alpha_{i\to j},\hat{x}_{i\to j}\}$
by passing the upstream messages $\{m_{k\to i}\}_{k\in\mathcal{N}_{i}\backslash j}$.
Here the message function $C_{i\to j}(x_{ij})$ is always smooth.
The messages $\alpha_{i\to j},\hat{x}_{i\to j}$ are computed as~\cite{Wong2007,Yeung2010,Rebeschini2019}
\begin{align}
\alpha_{i\to j} & =\frac{1}{\sum_{k\in\mathcal{N}_{i}\backslash j}\alpha_{k\to i}^{-1}}+r_{ij},\label{eq:message_alpha_ij_flow}\\
\hat{x}_{i\to j} & =\frac{\Lambda_{i}+\sum_{k\in\mathcal{N}_{i}\backslash j}\hat{x}_{k\to i}}{1+r_{ij}\sum_{k\in\mathcal{N}_{i}\backslash j}\alpha_{k\to i}^{-1}}.\label{eq:message_hat_xij_flow}
\end{align}

\subsubsection{The Reference Node $\mathcal{D}$\label{subsec:undirected_reference_node_D}}

Similar to the MP algorithm in routing games, there are two methods
to deal with possible boundary conditions of the reference node $\mathcal{D}$:
\begin{itemize}
\item Method I: Since node $\mathcal{D}$ has no constraints on its adjacent
flows, it will absorb all incoming flows resulting in $C_{\mathcal{D}\to j}(x_{\mathcal{D}j})=0$ and
\begin{equation}
\alpha_{\mathcal{D}\to j}=0,\quad\hat{x}_{\mathcal{D}\to j}=0.\label{eq:alpha_hatx_at_D_undirected}
\end{equation}
In this treatment, the Lagrange multiplier of node $\mathcal{D}$
can be set to $\mu_{\mathcal{D}}=0$ in the dual problem (as node
$\mathcal{D}$ is unconstrained), which corresponds to a grounded
node in the electric network interpretation of the problem.
\item Method II: Alternatively, one can set an explicit constraint on the
flows $\{x_{\mathcal{D}j}\}$ to the reference node $\mathcal{D}$
\begin{equation}
R_{\mathcal{D}}:=\Lambda_{\mathcal{D}}+\sum_{j\in\mathcal{N}_{\mathcal{D}}}x_{\mathcal{D}j}=0,
\end{equation}
where $\Lambda_{\mathcal{D}}=-\sum_{i\neq\mathcal{D}}\Lambda_{i}$.
Then the messages from the reference node $\mathcal{D}$ are calculated
in the same way as for other nodes.
\end{itemize}
Similar to the routing games, Method II results in an MP algorithm with
a faster convergence rate, but it may fail to converge for some graphs
while Method I can still provide valid solutions.

Method II is used in the experiments for undirected flow networks
in the main text.

\subsubsection{Computation of the Equilibrium Flows from Messages}

Upon convergence of the messages, the equilibrium flow $x_{ij}^{*}$
can be obtained by minimizing the edgewise full cost $C_{ij}^{\text{full}}(x_{ij})=C_{i\to j}(x_{ij})+C_{j\to i}(x_{ij})-\frac{1}{2}r_{ij}x_{ij}^{2}$,
giving rise to
\begin{equation}
x_{ij}^{*}=\frac{\alpha_{j\to i}\hat{x}_{j\to i}-\alpha_{i\to j}\hat{x}_{i\to j}}{\alpha_{i\to j}+\alpha_{j\to i}-r_{ij}}.\label{eq:undirected_equilibrium_flow_MP}
\end{equation}

\subsubsection{Results of the MP Algorithm on Undirected Flow Networks}

\begin{figure}
\includegraphics[scale=1.5]{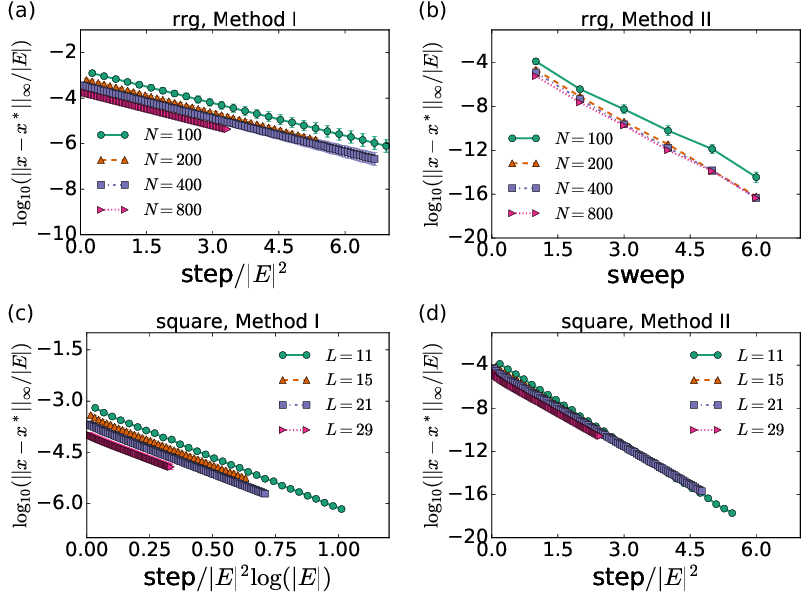}

\caption{MP algorithm in undirected flow networks converges to
the equilibrium flows $\boldsymbol{x}^{*}$. Random regular graphs
(RRG) of degree $3$ are considered in (a)(b), while square lattices
(size $L\times L$) are considered in (c)(d). Method I in Sec.~\ref{subsec:undirected_reference_node_D}
is used in (a)(c) to treat the reference node $\mathcal{D}$, while
Method II is used in (b)(d). In (b), each sweep comprises $40|E|$
local MP updates.\label{fig:undirected_flow_net_lower_MP}}

\end{figure}

In Fig.~\ref{fig:undirected_flow_net_lower_MP}, we demonstrate the
performance of the MP algorithms in undirected flow networks. The
MP algorithms converge in different networks, including square lattices
with many short loops. The iterations needed to obtained a given precision
seems to depend on the topologies of the networks, e.g., square lattices
appear to converge slower than random regular graphs. The method used to
treat the boundary reference node $\mathcal{D}$ also impacts on the
number of iterations needed, where it is observed that in general
Method II makes MP converge faster than Method I. We conjecture that
the influence of single-node boundary conditions in Eq.~(\ref{eq:alpha_hatx_at_D_undirected})
takes more iteration steps to diffuse messages to the bulk of the network. 

\subsection{MP For Bilevel Optimization}

The bilevel optimization problem on undirected networks aims to
tune the flows of targeted edges $\mathcal{T}$ such that they
exceed or drop below certain limits, depending on the application.
We consider the former case, where the goal is to control the flows on
$\mathcal{T}$ such that $\rho_{ij}(x_{ij})=\frac{|x_{ij}|-|x_{ij}^{0}|}{|x_{ij}^{0}|}-\theta\geq0,\forall(i,j)\in\mathcal{T}$
(with $x_{ij}^{0}$ being the flow before tuning). The task in the
upper-level is to minimize the objective
\begin{equation}
\min_{\boldsymbol{r}}\mathcal{O}(\boldsymbol{x}^{*}(\boldsymbol{r}))  = \sum_{(i,j)\in\mathcal{T}} \mathcal{O}_{ij}(x^{*}_{ij}(\boldsymbol{r}))  :=\sum_{(i,j)\in\mathcal{T}}-\rho_{ij}(x_{ij}^{*}(r))\Theta\big(-\rho_{ij}(x_{ij}^{*}(r))\big).
\end{equation}

As mentioned in the main text, the impact of the variation of the
control parameters $r_{ij}$ on the upper-level objective $\mathcal{O}$
is mediated through the messages $m_{i\to j}\in\{\alpha_{i\to j},\hat{x}_{i\to j}\}$
along the pathways from the targeted edges to edge $(i,j)$. For a
targeted edge $(p,q)\in\mathcal{T}$, the boundary conditions of the gradient with respect to the messages is given by $ \frac{\partial\mathcal{O}_{pq}}{\partial m_{p\to q}}  =\frac{\partial\mathcal{O}_{p q}}{\partial x_{pq}^{*}}\frac{\partial x_{pq}^{*}}{\partial m_{p\to q}}$,
where the components admit the following expressions 
\begin{align}
\frac{\partial\mathcal{O}_{pq}}{\partial x_{pq}^{*}}  & =  -\Theta\big(-\rho_{pq}(x_{pq}^{*})\big)\frac{\text{sgn}(x_{pq}^{*})}{|x_{pq}^{0}|},\\
\frac{\partial x_{pq}^{*}}{\partial\alpha_{p\to q}} & =\frac{-\hat{x}_{p\to q}}{\alpha_{p\to q}+\alpha_{q\to p}-r_{pq}}-\frac{-x_{pq}^{*}}{\alpha_{p\to q}+\alpha_{q\to p}-r_{pq}},\\
\frac{\partial x_{pq}^{*}}{\partial\hat{x}_{p\to q}} & =\frac{-\alpha_{p\to q}}{\alpha_{p\to q}+\alpha_{q\to p}-r_{pq}}.
\end{align}

The gradient with respect to the control parameter $r_{pq}$ on the targeted edge $(p, q)$ is computed as 
\begin{equation}
\frac{\partial\mathcal{O}_{pq}}{\partial r_{pq}} = \frac{\partial\mathcal{O}_{pq}}{\partial x_{pq}^{*}} \bigg[ \frac{\partial x_{pq}^{*}}{\partial r_{pq}} + \sum_{m\in\{\alpha,\hat{x}\}}\bigg( \frac{\partial x^{*}_{pq}}{\partial m_{p\to q}} \frac{\partial m_{p\to q}}{\partial r_{pq}} + \frac{\partial x^{*}_{pq}}{\partial m_{q\to p}} \frac{\partial m_{q\to p}}{\partial r_{pq}}\bigg) \bigg].
\label{eq:dO_dr_pq}
\end{equation}

For a non-targeted edge $(k,i)\notin\mathcal{T}$, we need to first evaluate $\frac{\partial\mathcal{O}}{\partial m_{k\to i}}$, which can be obtained by summing the gradients on its downstream edges $\{i\to l|l\in\mathcal{N}_{i}\backslash k\}$, computed as 
\begin{align}
\frac{\partial\mathcal{O}}{\partial m_{k\to i}} & = \sum_{l\in\mathcal{N}_{i}\backslash k}\,\,\sum_{m_{i\to l}\in\{\alpha_{i\to l},\hat{x}_{i\to l}\}}\frac{\partial\mathcal{O}}{\partial m_{i\to l}}\frac{\partial m_{i\to l}}{\partial m_{k\to i}} \nonumber \\
& = \sum_{(p,q) \in \mathcal{T}} \sum_{l\in\mathcal{N}_{i}\backslash k}\,\,\sum_{m_{i\to l}\in\{\alpha_{i\to l},\hat{x}_{i\to l}\}}\frac{\partial\mathcal{O}_{pq}}{\partial m_{i\to l}}\frac{\partial m_{i\to l}}{\partial m_{k\to i}}.
\end{align}
where the gradient propagation of messages $\frac{\partial m_{i\to l}}{\partial m_{k\to i}}$
admits the following forms
\begin{align}
\frac{\partial\alpha_{i\to l}}{\partial\alpha_{k\to i}} & =\frac{\alpha_{k\to i}^{-2}}{\big(\sum_{n\in\mathcal{N}_{i}\backslash l}\alpha_{n\to i}^{-1}\big)^{2}},\\
\frac{\partial\alpha_{i\to l}}{\partial\hat{x}_{k\to i}} & =0,\\
\frac{\partial\hat{x}_{i\to l}}{\partial\alpha_{k\to i}} & =\frac{\hat{x}_{i\to l}r_{il}\alpha_{k\to i}^{-2}}{1+r_{il}\sum_{n\in\mathcal{N}_{i}\backslash l}\alpha_{n\to i}^{-1}},\\
\frac{\partial\hat{x}_{i\to l}}{\partial\hat{x}_{k\to i}} & =\frac{1}{1+r_{il}\sum_{n\in\mathcal{N}_{i}\backslash l}\alpha_{n\to i}^{-1}}.
\end{align}

The resulting gradient w.r.t. to the
control parameter $r_{ki}$ (note that edge $(k,i)\notin\mathcal{T}$)
can be computed as
\begin{align}
& \frac{\partial\mathcal{O}}{\partial r_{ki}} = \sum_{(p,q)\in\mathcal{T}}  \frac{\partial\mathcal{O}_{pq}}{\partial r_{ki}} \ =  \sum_{(p,q)\in\mathcal{T}} \bigg[ \frac{\partial\mathcal{O}_{pq}}{\partial r_{k \to i}} + \frac{\partial\mathcal{O}_{pq}}{\partial r_{i \to k}} \bigg]  \\
& \frac{\partial\mathcal{O}_{pq}}{\partial r_{k \to i}}  =  \sum_{m\in\{\alpha,\hat{x}\}} \frac{\partial\mathcal{O}_{pq}}{\partial m_{k\to i}}\frac{\partial m_{k\to i}}{\partial r_{ki}},
\label{eq:dO_dr_ki}
\end{align}
where we have broken down the gradient as $\frac{\partial \mathcal{O}_{pq}}{\partial r_{ki}} = \frac{\partial \mathcal{O}_{pq}}{\partial r_{k \to i}} + \frac{\partial \mathcal{O}_{pq}}{\partial r_{i \to k}}$ for programming convenience.

In Eqs. (\ref{eq:dO_dr_pq}) and (\ref{eq:dO_dr_ki}), the terms $\frac{\partial m_{k\to i}}{\partial r_{ki}}$
are computed as
\begin{align}
\frac{\partial\alpha_{k\to i}}{\partial r_{ki}} & =1,\\
\frac{\partial\hat{x}_{k\to i}}{\partial r_{ki}} & =\hat{x}_{k\to i}\frac{-\sum_{n\in\mathcal{N}_{k}\backslash i}\alpha_{n\to k}^{-1}}{1+r_{ki}\sum_{n\in\mathcal{N}_{k}\backslash i}\alpha_{n\to k}^{-1}},
\end{align}
which closes the equations for the gradient computations.

In this formalism, we keep track of the influence of the control parameter $r_{ij}$ on a specific targeted edge $(p, q) \in \mathcal{T}$ in the form of $\frac{\partial \mathcal{O}_{pq}}{\partial r_{ij}}$.
The algorithmic complexity scales linearly with the number of targeted edges $|\mathcal{T}|$. This is acceptable as long as $|\mathcal{T}|$ is small. Further reducing the complexity (e.g., by considering the influence of the control parameter $r_{ij}$  on the overall objective function) is an interesting future research direction, which requires a more careful treatment of the boundary condition of the gradients.

We also consider the constraints on the control parameters, being in the range $r_{ij} \in [0.9, 1.1]$. More complex constraints can also be considered, e.g., by introducing a penalty function to the global objective, which will yield additional terms in the message-passing equations.

The resulting bilevel MP algorithms for flow control is summarized in Algorithm~\ref{alg:MP_flow_control_bilevel}.

\begin{algorithm}
\SetNoFillComment
\caption{Message-passing algorithm for (undirected) flow control (bilevel, Method II to treat the reference node $\mathcal{D}$ )}\label{alg:MP_flow_control_bilevel}

\KwIn{Undirected network $G(V,E)$ (pre-processed to remove dangling nodes), node parameters $\{ \Lambda_i \}$ defining the resources, maximal number of iterations $T$, time interval $t_{\text{update\_intv}}$ for performing gradient descent.}

\BlankLine
Initialize the messages $\{ \alpha_{i \to j}, \hat{x}_{i \to j} \}$ randomly.

\BlankLine
Initialize the gradients $\{ \partial \mathcal{O}_{pq} / \partial \alpha_{i \to j}, \partial \mathcal{O}_{pq} / \partial \hat{x}_{i \to j}, \partial \mathcal{O}_{pq} / \partial r_{i \to j}, \mathcal{O}_{pq} / \partial r_{pq} \}$ randomly.

\BlankLine
Initialize the edge control parameters $\{ r_{ij} \}$.

\BlankLine
\For{$t \textnormal{ in } 1:T$}{

\BlankLine
Randomly select a node $i$ and one of its adjacent node $j \in \mathcal{N}_i$.

\BlankLine
Update the messages $m_{i \to j} \in \{ \alpha_{i \to j}, \hat{x}_{i \to j} \}$ as:

\BlankLine
$\qquad \alpha_{i \to j} \leftarrow \frac{1}{\sum_{k \in \mathcal{N}_i \backslash j} \alpha_{k \to i}^{-1}}+r_{i j}$, \qquad $\hat{x}_{i \to j} \leftarrow \frac{\Lambda_i+\sum_{k \in \mathcal{N}_i \backslash j} \hat{x}_{k \to i}}{1+r_{i j} \sum_{k \in \mathcal{N}_i \backslash j} \alpha_{k \to i}^{-1}}$.

\BlankLine
\BlankLine

\For{ $(p,q) \in \mathcal{T}$ }{
\BlankLine

\uIf{$(i,j) = (p, q)$ \textnormal{or} $(j, i) = (p, q)$}{
\BlankLine
\tcc{$(i,j)$ is a targeted edge}
\BlankLine
Update the gradients $\partial \mathcal{O}_{pq} / \partial m_{p \to q}, \partial \mathcal{O}_{pq} / \partial m_{q \to p}, \partial \mathcal{O}_{pq} / \partial r_{pq}$ as

\BlankLine
$\qquad \frac{\partial \mathcal{O}_{p q}}{\partial m_{p \rightarrow q}} \leftarrow \frac{\partial \mathcal{O}_{p q}}{\partial x_{p q}^*} \frac{\partial x_{p q}^*}{\partial m_{p \rightarrow q}}$

\BlankLine
$\qquad \frac{\partial \mathcal{O}_{p q}}{\partial m_{q \rightarrow p}} \leftarrow \frac{\partial \mathcal{O}_{p q}}{\partial x_{p q}^*} \frac{\partial x_{p q}^*}{\partial m_{q \rightarrow p}}$

\BlankLine
$\qquad \frac{\partial \mathcal{O}_{p q}}{\partial r_{p q}} \leftarrow \frac{\partial \mathcal{O}_{p q}}{\partial x_{p q}^*}\left[\frac{\partial x_{p q}^*}{\partial r_{p q}}+\sum_{m \in\{\alpha, \hat{x}\}}\left(\frac{\partial x_{p q}^*}{\partial m_{p \rightarrow q}} \frac{\partial m_{p \rightarrow q}}{\partial r_{p q}}+\frac{\partial x_{p q}^*}{\partial m_{q \rightarrow p}} \frac{\partial m_{q \rightarrow p}}{\partial r_{p q}}\right)\right]$
\BlankLine
}

\Else{
\BlankLine
\tcc{$(i,j)$ is a non-targeted edge}
\BlankLine
Update the gradients $\partial \mathcal{O}_{pq} / \partial m_{i \to j}, \partial \mathcal{O}_{pq} / \partial r_{ij}, \forall (p,q) \in \mathcal{T}$ as

\BlankLine
$\qquad \frac{\partial \mathcal{O}_{pq}}{\partial m_{i \to j}} \leftarrow \sum_{l \in \mathcal{N}_j \backslash i} \sum_{m_{j \to l} \in\left\{\alpha_{j \to l}, \hat{x}_{j \to l}\right\}} \frac{\partial \mathcal{O}_{p q}}{\partial m_{j \to l}} \frac{\partial m_{j \to l}}{\partial m_{i \to j}}$,

\BlankLine
$\qquad \frac{\partial \mathcal{O}_{pq}}{\partial r_{i \to j}} \leftarrow \sum_{m_{i \to j} \in\left\{\alpha_{i \to j}, \hat{x}_{i \to j}\right\}} \frac{\partial \mathcal{O}_{p q}}{\partial m_{i \rightarrow j}} \frac{\partial m_{i \rightarrow j}}{\partial r_{i j}}$.
\BlankLine
}

} 

\BlankLine
\BlankLine
\If{$t \operatorname{mod} t_{\textnormal{update\_intv}} = 0$}{
\BlankLine
$\frac{\partial \mathcal{O}}{\partial r_{ij}} \leftarrow 0$. $\qquad$  \tcc{To compute the full gradient:}

\For{ $(p,q) \in \mathcal{T}$ }{
\BlankLine
\uIf{$(i,j) = (p, q)$ \textnormal{or} $(j, i) = (p, q)$}{
\BlankLine
$\frac{\partial \mathcal{O}}{\partial r_{ij}} \leftarrow \frac{\partial \mathcal{O}}{\partial r_{ij}} + \frac{\partial \mathcal{O}_{pq}}{\partial r_{pq}}$.
}

\Else{
\BlankLine
$\frac{\partial \mathcal{O}}{\partial r_{ij}} \leftarrow \frac{\partial \mathcal{O}}{\partial r_{ij}} + \frac{\partial \mathcal{O}_{pq}}{\partial r_{i \to j}} + \frac{\partial \mathcal{O}_{pq}}{\partial r_{j \to i}}$.
}

} 

\BlankLine
$r_{ij} \leftarrow r_{ij} - s \frac{\partial \mathcal{O}}{\partial r_{ij}} \qquad$ \tcc{Update the control parameter with step size $s$.}

$r_{ij} \leftarrow \min( \max(r_{ij}, r_{\min}), r_{\max} ) \qquad$ \tcc{To respect the box constraints of the control parameters.}
\BlankLine
}

\BlankLine
\BlankLine
\If{\textnormal{messages converge}}{
Exit the for loop.
}

} 

\BlankLine

For each edge $e$, compute the equilibrium flow as $x_{ij}^*=\frac{\alpha_{j \rightarrow i} \hat{x}_{j \rightarrow i}-\alpha_{i \rightarrow j} \hat{x}_{i \rightarrow j}}{\alpha_{i \rightarrow j}+\alpha_{j \rightarrow i}-r_{i j}}$.

\BlankLine
\KwOut{Convergence status, equilibrium flows $\{ x^{*}_{ij} \}$ and optimal control parameters $\{ r_{ij} \}$}
\end{algorithm}

\subsection{The Global Optimization Approach}

In this section, we provide details of the global optimization approach
on undirected flow networks used in the main text. As mentioned before, e.g., in Sec.~\ref{sec:MPLowLevelOpt},
we have $x_{ij}^{*}=\frac{\mu_{j}^{*}-\mu_{i}^{*}}{r_{ij}}$, where
$\boldsymbol{\mu}^{*}=L(\boldsymbol{r})^{\dagger}\boldsymbol{\Lambda}$
and $L(\boldsymbol{r})^{\dagger}$ is the pseudo-inverse of the Laplacian
matrix $L(\boldsymbol{r})$ defined in Eq.~(\ref{eq:lapl_matrix_def}).
To compute the gradient $\frac{\partial\mathcal{O}}{\partial r_{ij}}$,
we need to evaluate the response of $\boldsymbol{\mu}^{*}$ to the
variation of the control parameters $\boldsymbol{r}$, i.e.,
\begin{equation}
\frac{\partial\mathcal{O}}{\partial r_{ij}}=\sum_{(p,q)\in\mathcal{T}}\frac{\partial\mathcal{O}}{\partial x_{pq}^{*}}\bigg[\frac{1}{r_{pq}}\big(\frac{\partial\mu_{q}^{*}}{\partial r_{ij}}-\frac{\partial\mu_{p}^{*}}{\partial r_{ij}}\big)-\frac{1}{r_{pq}^{2}}(\mu_{q}^{*}-\mu_{p}^{*})\delta_{(i,j),(p,q)}\bigg].
\end{equation}

Furthermore, it requires to compute $\frac{\partial L(\boldsymbol{r})^{\dagger}}{\partial r_{ij}}$
for $\frac{\partial\boldsymbol{\mu}^{*}}{\partial r_{ij}}$. Assuming
the underlying graph will not fragment into multiple components
when adapting $\boldsymbol{r}$, $L(\boldsymbol{r})$ has a constant
rank, then we have~\cite{Golub1973}
\begin{align}
\frac{\partial L(\boldsymbol{r})^{\dagger}}{\partial r_{ij}}= & -L^{\dagger}\frac{\partial L}{\partial r_{ij}}L^{\dagger}+L^{\dagger}(L^{\dagger})^{\top}\big(\frac{\partial L^{\top}}{\partial r_{ij}}\big)(I-LL^{\dagger})\nonumber \\
 & +(I-L^{\dagger}L)\big(\frac{\partial L^{\top}}{\partial r_{ij}}\big)(L^{\dagger})^{\top}L^{\dagger}.
\end{align}
Using the property of the pseudo-inverse of the Laplacian $L^{\dagger}L=LL^{\dagger}=I-\frac{1}{N}\boldsymbol{1}\boldsymbol{1}^{\top}$
and $\frac{\partial L}{\partial r_{ij}}\cdot\boldsymbol{1}=0$ (with
$\boldsymbol{1}$ as the all-one vector), we have
\begin{equation}
\frac{\partial L(\boldsymbol{r})^{\dagger}}{\partial r_{ij}}=-L^{\dagger}\frac{\partial L}{\partial r_{ij}}L^{\dagger},
\end{equation}
which closes the equations for calculating the gradients.

\clearpage

\def\bibsection{\section*{\refname}} 
\bibliography{reference}